\documentclass[12pt]{article}
\usepackage{graphicx}
\usepackage{epsfig}
\usepackage{amssymb, amsmath}
\usepackage{amsthm}
\usepackage{setspace}
\usepackage[top = 1in, bottom = 1in, left = 1.2in, right = 1.2in]{geometry}
\newtheorem{theorem}{Theorem}[section]
\newtheorem{corollary}[theorem]{Corollary}
\newtheorem{lemma}[theorem]{Lemma}
\newtheorem{prop}[theorem]{Proposition}
\newtheorem{remark}[theorem]{Remark}

\newenvironment{subproof}{\noindent{\em Pf:  }}

\newcommand {\aplus}[1] {\ensuremath{a_{#1}^{+}}}
\newcommand {\amin}[1] {\ensuremath{a_{#1}^{-}}}

\newcommand {\aplusk}[1] {\ensuremath{a_{#1}^{+k}}}
\newcommand {\amink}[1] {\ensuremath{a_{#1}^{-k}}}

\newcommand{\ddt}{\ensuremath{\frac{\partial}{\partial t}}}
\newcommand{\dt}[1]{\ensuremath{\partial_{t} #1}}

\newcommand {\K}{\ensuremath{\mathcal{K}}}
\newcommand {\lambdah}[1] {\ensuremath{\lambda_{#1}^{h}}}
\newcommand{\Laph}{\ensuremath{\Delta^{h}}}
\newcommand{\letilde}{\ensuremath{\lesssim}}
\newcommand{\Lh}[1]{\ensuremath{L_{h}^{#1}}}
\newcommand{\mN}{\ensuremath{\mathcal{N}}}

\newcommand{\Nbold}[1]{\ensuremath{\mathbb{N}^{#1}}}
\newcommand{\fracpartial}[2]{\ensuremath{\frac{\partial {#1}}{\partial {#2} }}}
\newcommand{\phih}{\ensuremath{\phi^{h}}}
\newcommand{\phis}{\ensuremath{\overline{\phi}}}
\newcommand{\PmN}{\ensuremath{\mathbf{P}_{T_{u_{j}}\mathcal{N}}}}

\newcommand {\qplus}[1] {\ensuremath{q_{#1}^{+}}}
\newcommand {\qmin}[1] {\ensuremath{q_{#1}^{-}}}

\newcommand {\qplusa}[1] {\ensuremath{q_{#1}^{+\alpha}}}
\newcommand {\qmina}[1] {\ensuremath{q_{#1}^{-\alpha}}}
\newcommand {\qplusk}[1] {\ensuremath{q_{#1}^{+k}}}
\newcommand {\qmink}[1] {\ensuremath{q_{#1}^{-k}}}
\newcommand {\qpluski}[2] {\ensuremath{q_{#1}^{+k, #2}}}
\newcommand {\qminki}[2] {\ensuremath{q_{#1}^{-k, #2}}}
\newcommand{\R}[1]{\ensuremath{\mathbb{R}^{#1}}}
\newcommand{\Rplus}{\ensuremath{\mathbb{R}_{+}}}

\newcommand{\TN}{\ensuremath{T_{u_{j}}\mathcal{N}}}

\newcommand{\TpN}{\ensuremath{T_{p}\mathcal{N}}}
\newcommand{\TyN}{\ensuremath{T_{y}\mathcal{N}}}
\newcommand {\U}{\ensuremath{\mathcal{U}}}
\newcommand {\uh}[1] {\ensuremath{u_{#1}^{h}}}

\newcommand {\vplus}[1] {\ensuremath{v_{#1}^{+}}}
\newcommand {\vmin}[1] {\ensuremath{v_{#1}^{-}}}
\newcommand {\vplusa}[1] {\ensuremath{v_{#1}^{+i}}}
\newcommand {\vmina}[1] {\ensuremath{v_{#1}^{-i}}}
\newcommand {\vplusk}[1] {\ensuremath{v_{#1}^{+k}}}
\newcommand {\vmink}[1] {\ensuremath{v_{#1}^{-k}}}
\newcommand {\wplus}[1] {\ensuremath{w_{#1}^{+}}}
\newcommand {\wmin}[1] {\ensuremath{w_{#1}^{-}}}

\newcommand {\wplusk}[1] {\ensuremath{w_{#1}^{+k}}}
\newcommand {\wmink}[1] {\ensuremath{w_{#1}^{-k}}}

\begin{document}
\title{\LARGE The Construction of a Partially Regular Solution to the Landau-Lifshitz-Gilbert Equation in \R{2}}
\author{Joy Ko \\ \it Brown University }
\singlespacing
\maketitle
\numberwithin{equation}{section}

\begin{abstract}
We establish a framework to construct a global solution in the space of finite energy to a general form of the Landau-Lifshitz-Gilbert equation in \R{2}.   Our characterization yields a partially regular solution, smooth away   from a $2$-dimensional locally finite Hausdorff measure set.  This construction relies on approximation by discretization, using the special geometry to express an equivalent system whose highest order terms are linear and  the translation of  the  machinery of linear estimates on the fundamental solution from the continuous setting into the discrete setting.  This method is quite general and accommodates more general geometries involving targets that are compact smooth hypersurfaces.  
\end{abstract}

\section{Introduction}

Micromagnetics, a model based on the work of Landau and Lifshitz, offers a description of magnetic behavior in ferromagnetic materials.   This continuum model is based on finding critical points of the micromagnetic energy associated to the magnetic moment represented by the field  $u: \Omega \rightarrow S^2$:
\[
E[u] = \int_{\Omega} |\nabla u|^{2} + \kappa\int_{\Omega}\phi(u) + \int_{{\mathbb{R}}^{3}}|\nabla M|^{2} dx - 2\int_{\Omega} h_{\mathrm ext}\cdot u
\]
where $\Omega$ is the region occupied by the ferromagnet and  all the physical constants have been normalized to be $1$.  These four terms are known as the exchange, anisotropic, magnetostatic, and  external energies respectively. 
The exchange energy penalizes spatial variation 
over macroscopic distances. The anisotropic term favors special directions of the magnetization. 
The magnetostatic energy represents the work required to build up the macroscopic body
by bringing its magnetic moments from infinity to their eventual position;  $M$ is nonlocal
and is  defined by $\Delta M = \mbox{div } u$ where $u$ is extended to be $0$ outside
of $\Omega$ and the equation is to be understood in the sense of distributions.
Finally, the external energy favors magnetization aligned with an external applied field.
The first variation of the energy with respect to $u$ is called the effective field:
\[
h_{eff} = -\frac{\delta E[u]}{\delta u}.
\]
Constrained to the set of candidate solutions $\{u: \Omega  \rightarrow S^{2} \}$, 
the resulting dynamic equations associated to this energy are, in dimensionless form,  given by
\[
\partial_{t}u =  u \wedge h_{\mathrm eff}, 
\]
where $\wedge$ denotes the wedge product in $\R{3}$.
To incorporate the Gilbert damping law, whose origin lies in the observation that such systems reach
equilibrium and must have decreasing energy over time, a dissipative term can be added on resulting in the 
following equation:
\[
\partial_{t}u = u \wedge h_{\mathrm eff} -  \alpha u \wedge (u \wedge h_{\mathrm eff}).  
\]
In the physical setting, the magnetization configuration is
the compromise that results from the competition of all of these terms to minimize the energy.  However, 
many relevant features of solutions to the full 
equation are captured by retaining the term associated to the exchange energy in which case
 $h_{\mathrm eff} = \Delta u$.    For the question we address here, in particular, the other terms can be considered lower order  terms.  
An analytically interesting variant of the problem results from attempting to enlarge the class of  target spaces beyond  $S^{2}$ to hypersurfaces with unit normal $\nu$.  
Although this variant does not hold common currency in micromagnetics, this might be a reasonable means of incorporating anisotropic contributions. 

The Cauchy problem for this generalized Landau-Lifshitz-Gilbert equation (LLG)  is the problem of finding $u$, given initial data $f : \R{2} \rightarrow \mathcal{N}$ for $\mathcal{N}$ a compact smooth hypersurface in $\R{3}$ with normal vector $\nu$,  satisfying
\begin{equation}
\label{llgequation}
\left\{\begin{array}{lcl}
\displaystyle
\partial_{t}u = \nu(u) \wedge \Delta u - \alpha \nu(u) \wedge (\nu(u) \wedge \Delta u) \\
u(x, 0) = f
\end{array}\right., 
\end{equation}
where $\alpha > 0$.  
We will refer to the first term as the Schr\"odinger term (this terminology will be clarified later) and the second as the damping term.    When only the damping term is present, this equation is the \emph{harmonic map heat flow problem}, which will later be a focus of discussion.  Although standard usage of LLG refers to (\ref{llgequation}) in the special case $\mathcal{N} = S^2$,  we adopt  here a usage to accommodate this general class of target spaces.

There are several standard  forms of ($\ref{llgequation}$) which are equivalent for smooth solutions.  The first results from the vector identity $-\xi \wedge (v \wedge \xi) = -v + (v, \xi) v$
which holds for $v$ a unit vector.  From this, 
($\ref{llgequation}$) can be written as
\begin{equation}
\label{llgform2}
\partial_{t}u = \nu(u) \wedge \Delta u + \alpha (\Delta u + (\nabla \nu(u) \cdot \nabla u) \nu(u))
\end{equation}
From (\ref{llgequation}) and (\ref{llgform2}), two additional formulations follow:
\begin{equation}
\label{llgform3}
\alpha \partial_{t}u -  \nu(u) \wedge \partial_t u = (1 + \alpha^2) (\Delta u + (\nabla \nu(u)\cdot \nabla u) \nu(u))
\end{equation}
\begin{equation}
\label{llgform4}
\partial_{t}u + \alpha \nu(u) \wedge \partial_t u = (1 + \alpha^2) \nu(u) \wedge \Delta u
\end{equation}

 The question that we address in this paper is whether global solutions exist in the class of finite Dirichlet energy.  
The main result in this paper is the following:

\begin{theorem}
\label{fellgthm1}
For any $f \in H^{1}(\R{2}, \mathcal{N})$, there exists, for $t > 0$, a solution to (\ref{llgequation}) which is smooth away from a singular set that has locally finite $2$-dimensonal Hausdorff measure with respect to the parabolic metric.
\end{theorem}

Weak solutions, even partially regular ones, have been shown to exist for such equations in the case $\mathcal{N} = S^2$.   The techniques developed, however, often exploit the special geometry of $S^2$ and side-step the natural difficulties that the general problem poses. Consequently, these techniques are not easily adaptable to the  setting that we consider here.  One category of results showing existence of weak solutions has made crucial use of the fact that  LLG can be written in divergence form in the case $\mathcal{N} = S^2$;  this in turn gives rise to a definition of weak solution.  Amongst these is the work of Alouges and Soyeur  \cite{alougessoyeur} who have shown that in three dimensions,  energy bounds are sufficient for existence of such a weak solution.   Guo and Hong \cite{guohong} successfully carried through the argument that Struwe in \cite{struwe} employed for the harmonic map heat flow to exhibit a \emph{Struwe solution}, a partially regular solution that satisfies an energy inequality and is smooth away from a finite set of point singularities.   In this dimension, Struwe shows that the existence of smooth solutions in addition to energy estimates is 
sufficient for the characterization of such a partially regular solution.  In the case of the sphere, 
the problematic term $u \wedge \Delta u$ does not preclude energy estimates since multiplication by $u$ and $\Delta u$ eliminates this term altogether.  In addition, they show that  the linear term in LLG is strictly parabolic and have used parabolic $L^{p}$ theory to establish local existence.  It is worthwhile to note that uniqueness in the class of partially regular solutions is open;  our partially regular solution may
be different from the Struwe solution.
Through private communication, we have learned that Melcher in \cite{melcher} has recently attained a characterization of partially regular solution similar to ours in three dimensions for $\mathcal{N} = S^2$ using a Ginzburg-Landau approximation.  

In the context of existing results, the main difficulties are already apparent.  
The nonlinearity in the highest order $\nu(u) \wedge \Delta u$ term poses the most apparent obstacle.  At first glance, it renders unavailable the machinery of partial regularity theory established for semilinear equations:   linear estimates on the fundamental solution  (e.g., $L^p-L^q$, Strichartz-type) and regularity arguments based on the inversion of the fundamental solution.    An additional problem that this term poses for general target $\mathcal{N}$  is that the local existence of smooth solutions is not just an immediate consequence of parabolic theory. 
 Crucial tools that are lost  include the maximum principle as well as a Bochner identity.   Lastly,  the divergence structure of the equation is lost when we depart from the sphere.
 
 Our method is based on approximating systems on a uniform spatial grid.  The idea is to show that  a sequence of solutions to this approximation converges to a solution to the original equation as the grid-size approaches zero.  The idea of approximating systems via discretization is not new;  in fact, one of the earliest works investigating ($\ref{llgequation}$) for the case $\alpha = 0$, $\mathcal{N} = S^2$ \cite{ssb} uses discrete approximations to show local existence.  However, extendings its scope to partial regularity is new.    While the techniques developed here are demonstrated for the specific setting of the LLG equations,  the generality of these techniques should be emphasized.  This key ingredients used are:
 (i) the construction of a suitable discretization of the system; 
(ii) establishing a suitable linearization of the discrete system; 
(iii) the construction of the discrete fundamental solution of the resulting linear operator and deriving appropriate linear mixed space-time estimates.  
 The setting that enables the discrete  constructions described by (i) and (iii) are the subject of Section $\ref{discretesetting}$.   

The advantages of this method address the difficulties of the problem.  By discretizing the system, existence of solutions to the resulting ODE is immediate once we attain uniform bounds;  there is no 
need to appeal to heavy machinery to conclude that smooth solutions exist.   The  transformation which results in a linear highest order term opens the door to the arsenal of linear
estimates and methods. The explicit construction of the discrete fundamental solution and the
resulting linear estimates permit higher derivative bounds that has no dependence on the special
structure of the sphere. 

Before employing the proposed method to construct a partially regular solutions for LLG,
 we illustrate in Section $\ref{hhp}$ this method for
the harmonic map heat flow problem, which involves only the damping term.  The same characterization of a partially regular solution has been achieved by other methods $\cite{cs}$ , but we consider this excursion to be worthwhile since it is a setting in which we can isolate the difficulties  that just one term poses  (makes use of ingredients (i) and (iii) above) before adding on the separate difficulty involving the
Schr\"odinger term.   In Section $\ref{llg}$, our main result is proved. 

\section{The Discrete Setting}
\label{discretesetting}
\subsection{Notation and Definitions}

Consider a uniform grid on $\R{d}$ determined by the lines $x_{i} = j_{i}h$,
$i = 1,2, \ldots, d$ where $j_{i}$ are integers. For vector valued functions $u^{h}$
defined on the grid,  the definition of scalar product and discrete $L^{p}$ norms that will be used  are:
\[
(u^{h}, v^{h})_{\Lh{2}} = 
   h^{d}\sum_{j} u^{h}_{j} v^{h}_{j}
\]
\[
\|u^{h}\|_{\Lh{p}} = (h^{d}\sum_{j} |u^{h}_{j}|^{p})^{1/p}, 
\]
where 
$u_{j}^{h} = u^{h}_{j_{1}, \ldots, j_{d}} = u(j_{1}h, \ldots, j_{d}h)$  and  
$u^{h}_{j_{i}+ 1} = u(j_{1}, \ldots, j_{i}h + h, \ldots,  j_{d}h)$. 
The  basic difference operations are
\[
D_{+i}u^{h}_{j} = \frac{u^{h}_{j_{i} + 1} -u^{h}_{j}}{h},  \ \
D_{-i}u^{h}_{j} = \frac{u^{h}_{j} -u^{h}_{j_{i}-1}}{h},  \ \
D_{0i}u^{h}_{j} = \frac{u^{h}_{j_{i} + 1} -u^{h}_{j_{i}-1}}{2h}, \ \
\]
and 
\[
\delta_{+i}^{2}u^{h}_{j} = D_{+i}D_{-i}u^{h}_{j} =  
D_{-i}D_{+i}u^{h}_{j} = \frac{1}{h}(D_{+i}-D_{-i})u^{h} 
 = \frac{u^{h}_{j_{i} + 1}  - 2u^{h}_{j} + u^{h}_{j_{i}-1}}{h^{2}}. 
\]
The discrete Laplacian $\Delta^{h}$ can be defined as 
\[
\Delta^{h} = \sum_{i = 1}^{d} \delta_{+i}^{2}.
\]
For some multi-index  $\alpha = (\alpha_{1}^{+}, \alpha_{1}^{-}, \ldots,\alpha_{d}^{+}, \alpha_{d}^{-})$, $\alpha_{k}^{+}, \alpha_k^- \in \Nbold{}$, 
we will use the notation
\[
|\alpha| = \sum_{i = 1}^{d} (\alpha_{i}^{+} +  \alpha_{i}^{-})
\]
and
\[
D^{\alpha} = \prod_{i = 1}^{d}(D_{+i})^{\alpha_{i}^{+}}(D_{-i})^{\alpha_{i}^{-}}.
\]
The discrete Sobolev spaces are defined as
\[
W_{h}^{k, p} \stackrel {\Delta}{=}  \{u^{h} |  \|u^{h}\|_{W_{h}^{k,p}} = \sum_{|\alpha| \leq k} \|D^{\alpha}u^{h}\|_{\Lh{p}} < \infty\}
\]
with $H_{h}^{k} = W_{h}^{k, 2}$.  Their  homogeneous counterparts are defined by
\[ 
\dot{W_{h}}^{k, p} \stackrel {\Delta}{=}  \{u^{h} |  \|u^{h}\|_{W_{h}^{k,p}} = \sum_{|\alpha| = k} \|D^{\alpha}u^{h}\|_{\Lh{p}} < \infty\}
\] 
with $\dot{H_{h}}^{k} = \dot{W_{h}}^{k, 2}$.

For integer $k$, we will use the notation $D^k u_j$ to denote the vector consisting of all difference derivatives $D^{\alpha} \uh{j}, |\alpha| = k$.   For the case $k = 1$, $D^{+}$  and $D^{-}$ will denote the 
approximate gradients, the vectors with components $D_{+k} \uh{j}$ and $D_{-k}\uh{j}$ respectively. 
\subsubsection*{Notable Imperfections of Discretization}
\label{defects}
Discretization introduces some imperfection which are important to highlight.   
The common operations of product differentiation and
integration by parts have a slight shift.  The forms of product differentiation that will be used are:
\begin{eqnarray*}
D_{+i}(u^{h}v^{h})_{j} 
& = &  u^{h}_{j} D_{+i}v^{h}_{j} + v^{h}_{j_{i} + 1} D_{+i}u^{h}_{j} \\
& = & v^{h}_{j} D_{+i}u^{h}_{j} + u^{h}_{j_{i} + 1} D_{+i}v^{h}_{j} \\
& = & \frac{1}{2}\{(u^{h}_{j_{i} + 1} + u^{h}_{j}D_{+i}v^{h}_{j}
  + (v^{h}_{j_{i} + 1} + v^{h}_{j})D_{+i}u^{h}_{j}\}.
\end{eqnarray*}
For functions $u^{h}, v^{h}$ both belonging to $\Lh{2}$, the summation by 
parts formula-- the discrete analog of integration by parts-- is given by
\[
(u^{h}, D_{+i}v^{h})_{\Lh{2}} =-(D_{-i}u^{h}, v^{h})_{\Lh{2}}.
\]
Another important imperfection that arises is that difference derivatives lose their `tangency'.  Consider
a map  $u$ from $\R{2}$ to $\mathcal{N}$, where $\mathcal{N}$ is a hypersurface in $\R{3}$ with unit normal $\nu$.  
Difference derivatives of $u$ no longer lie in the tangent plane $\TN$.   To negotiate some of the difficulty that this defect creates, we will often separate the portion of the difference derivative that lies in $\TN$ and that which is orthogonal.  

\subsection{Basic Inequalities and Some Useful Tools}
Most of the inequalities that will be used have long been established.   Amongst these are Holder's, Minkowski's and the $L^{p}$ interpolation inequalities which translate perfectly to the discrete setting.  
Discrete versions of the Sobolev embedding theorems were established by Ladyzhenskaya \cite{lady} using interpolation
operators.  In particular, we will make use of the following discrete Sobolev-interpolation inequality:
\begin{equation}
\label{sob-int}
\|D^{j}u^{h}\|_{\Lh{r}} \leq  C\|u^{h}\|_{\Lh{p}}^{1-\theta} \|D^{k}u^{h}\|_{\Lh{q}}^{\theta}, \mbox{     } 
   \frac{1}{r} -\frac{j}{n}= \frac{1-\theta}{p} + \theta(\frac{1}{q}-\frac{k}{n}).
\end{equation}
The interpolation operators that supply the proof for the above inequalities are useful and we mention
one which will play an immediate and a later role.  Confining our attention to the case $d = 2$ for the remainder of this section,
associate with $u$ (which is only defined on the uniform grid in the plane) the unique interpolated polynomial 
$p_{h}$ which matches the value of $u$ at each grid point and is of the
form $a_{0} + a_{1}x + a_{2}y + a_{3}xy$ in each square. Specifically,  for $j$ indexing pairs $(j_{1}, j_{2})$, 
\[
p_{h} = \sum_{j} (u_{j} + D_{+1}u_{j}(x-j_{1} h) + D_{+2}u_{j}(y-j_{2}h) + D_{+1}D_{+2}u_{j}(x-j_{1} h)
(y-j_{2}h))\chi_{\Box_{j}}(x,y),
\]
where $\chi_{\Box_{j}}(x,y)$ denotes the characteristic function on the square $K_{j}$.
Across any one side of a square, $p_{h}$ 
is linear, so we have that $p_{h}$ is continuous across interfaces.  $p_{h}$ is an $H^{1}$-interpolant
associated to $u$.    Of particular importance are the equivalence of the norms $\|p_{h}\|_{L^p}$ and $\|u^h\|_{\Lh{p}}$ as well as $\|\nabla p_{h}\|_{L^p}$ and $\|D^1 u^h\|_{\Lh{p}}$ for any $p \geq 1$.
Equipped with this interpolant, we can prove the following propositions.  The first is a localized Sobolev-interpolation inequality.  
\begin{prop}
\label{locsob-int}
Let $\Omega \subset \R{2}$ and  $\zeta \in C_{0}^{\infty}(\Omega)$.  Then
\[
\||u^{h}|^{2}\zeta\|_{\Lh{2s}(\R{2})} \leq C\|u^{h}\|_{\Lh{p}(\Omega)} \|D^{1}(u^{h}\zeta)\|_{\Lh{q}(\R{2})},  \mbox{       }   \frac{s+1}{2s} = \frac{1}{p} + \frac{1}{q}.
\]
\end{prop}
\begin{remark}
Without the $\zeta$, the quantity of interest is $\|u^{h}\|_{\Lh{4s}(\R{2})}^{2}$ and the proposed inequality
is given by $(\ref{sob-int})$ for $n = 2, r = 4s, j = 0, k = 1, \theta = \frac{1}{2}.$
\end{remark}
\begin{proof}
This inequality follows from the equivalence of the norms  $\|p_{h}\|_{L^p}$ and $\|u^h\|_{\Lh{p}}$ as well as $\|\nabla p_{h}\|_{L^p}$ and $\|D^1 u^h\|_{\Lh{p}}$ for any $p \geq 1$ and  the equivalent statement in the continuous setting, as proved in Appendix A. 
\end{proof}
\begin{prop}
\label{phlemma}
If $\|u^{h}\|_{H_{h}^{2}(\R{2})}< C$ independent of $h$, then there is a subsequence $\{p_{h_{k}}\}$ that converges strongly in $H^{1}(\R{2})$. 
\end{prop}
\begin{proof}
Notice that the boundedness of $p_{h}$ in $H^{1}$ immediately gives weak convergence in $H^{1}$.
To show strong $H^{1}$ convergence, we employ Ascoli-Arzela.  In particular, we must show that
$\{\nabla p_{h}\}$ are bounded and equicontinuous in $L^{2}$. 

We have $L^{2}$ bounds on $\nabla p_{h}$ by construction. To show equicontinuity, we must show that $\forall \epsilon > 0, \exists
\delta  = (\delta_{1}, \delta_{2}) = (k_{1}h, k_{2}h)$ such that 
$\|\nabla p_{h}(\cdot + \delta_{1}, \cdot + \delta_{2}) - \nabla p_{h}\|_{L^{2}} < c|\delta|.$
We have
\begin{eqnarray*}
\lefteqn{\|\nabla p_{h}(\cdot + \delta_{1}, \cdot + \delta_{2}) - \nabla p_{h}\|_{L^{2}}^{2}}\\
& = & \sum_{j_{1} ,j_{2}} \int_{\Box_{j_{1} ,j_{2}}}|(D_{+1}u_{j_{1} +k_{1},j_{2}}-D_{+1}u_{j_{1} ,j_{2}}) + (D_{+1}D_{+2}u_{j_{1} +k_{1},j_{2}}- \\ & & \mbox{} D_{+1}D_{+2}u_{j_{1} ,j_{2}})(y-j_{2}h)|^{2} +
|(D_{+2}u_{j_{1} ,j_{2}+k_{2}}-D_{+2}u_{j_{1} ,j_{2}}) + (D_{+1}D_{+2}u_{j_{1} ,j_{2}+k_{2}} - \\ & & \mbox{}-D_{+1}D_{+2}u_{j_{1} ,j_{2}})(x-j_{1}h)|^{2} dx dy\\
& \leq  & 2\sum_{j_{1} ,j_{2}} \int_{\Box_{j_{1} ,j_{2}}} |D_{+1}u_{j_{1} +k_{1},j_{2}}-D_{+1}u_{j_{1} ,j_{2}}|^{2} + |D_{+2}u_{j_{1} ,j_{2}+k_{2}}-D_{+2}u_{j_{1} ,j_{2}}|^{2} +\\ & & 
\mbox{} |D_{+1}D_{+2}u_{j_{1} +k_{1},j_{2}}-D_{+1}D_{+2}u_{j_{1} ,j_{2}}(y-j_{2}h)|^{2} +|D_{+1}D_{+2}u_{j_{1} ,j_{2}+k_{2}}- \\ & &
\mbox{} D_{+1}D_{+2}u_{j_{1} ,j_{2}}(x-j_{1}h)|^{2} dx dy\\
& = & 2(I + II).  
\end{eqnarray*}
To bound $I$, we can write the difference $D_{+1}u_{j_{1} +k_{1},j_{2}}-D_{+1}u_{j_{1} ,j_{2}}$ as
\[
D_{+1}u_{j_{1} +k_{1},j_{2}}-D_{+1}u_{j_{1} ,j_{2}} = \sum_{k = 1}^{k_{1}} (D_{+1}u_{j_{1} +k,j_{2}}-D_{+1}u_{j_{1} +k-1,j_{2}}).
\]
Since $(a_{1} + \cdots +a_{n})^{2} \leq 2n(a_{1}^{2} + \cdots  + a_{n}^{2})$, we have the bound
\[
|D_{+1}u_{j_{1} +k_{1},j_{2}}-D_{+1}u_{j_{1} ,j_{2}}|^{2} \leq 2k_{1}\sum_{k = 1}^{k_{1}} |D_{+1}u_{j_{1} +k,j_{2}}-D_{+1}u_{j_{1} +k-1,j_{2}}|^{2}.
\]
Similarly, we have the  bound
\[
|D_{+2}u_{j_{1} ,j_{2}+ k_{2}}-D_{+1}u_{j_{1} ,j_{2}}|^{2} \leq 2k_{2}\sum_{k = 1}^{k_{2}} 
|D_{+1}u_{j_{1} ,j_{2}+k}-D_{+1}u_{j_{1} ,j_{2}+k-1}|^{2}.
\]
Using the given information that $\|u^h\|_{H_h^2} < C$, we have in particular that
\[
\|\Delta u\|_{L_{h}^{2}}^{2} = \displaystyle \sum_{j_{1} ,j_{2}} |D_{+1}u_{j_{1} ,j_{2}}-D_{-1}u_{j_{1} ,j_{2}}|^{2} + 
|D_{+2}u_{j_{1} ,j_{2}}^{h}-D_{-2}u_{j_{1} ,j_{2}}^{h}|^{2} < C,
\]
which allows us to conclude that 
\begin{eqnarray*}
I & =  & \sum_{j_{1} ,j_{2}} \int_{\Box_{j_{1} ,j_{2}}} |D_{+1}u_{j_{1} +k_{1},j_{2}}-D_{+1}u_{j_{1} ,j_{2}}|^{2} + |D_{+2}u_{j_{1} ,j_{2}+k_{2}}-D_{+2}u_{j_{1} ,j_{2}}|^{2} \\
  & \leq & h^{2}\sum_{j_{1} ,j_{2}} |D_{+1}u_{j_{1} +k_{1},j_{2}}-D_{+1}u_{j_{1} ,j_{2}}|^{2} + |D_{+2}u_{j_{1} ,j_{2}+k_{2}}-D_{+2}u_{j_{1} ,j_{2}}|^{2} \\
  & \leq & h^{2}[2(k_{1}-1)\sum_{j_{1} ,j_{2}}\sum_{k = 1}^{k_{1}} |D_{+1}u_{j_{1} +k,j_{2}}-D_{+1}u_{j_{1} +k-1,j_{2}}|^{2} + \\ & & \mbox{} 
                 2(k_{2}-1)\sum_{j_{1} ,j_{2}}\sum_{k = 1}^{k_{2}} |D_{+1}u_{j_{1} ,j_{2}+k}-D_{+1}u_{j_{1} ,j_{2}+k-1}|^{2}] \\
  & \leq & 2h^{2}[(k_{1}-1)k_{1}C + (k_{2}-1)k_{2}C] \\
  & \leq & 2Ch^{2}(k_{1}^{2} + k_{2}^{2}) \\
  & \leq & c|\delta|^{2}.
\end{eqnarray*}
By assumption, we also have that $\sum_{j_{1} ,j_{2}}h^{2}  |D_{+1}D_{2} u_{j_{1} ,j_{2}}|^{2} \leq C$
so
\begin{eqnarray*}
II  & \leq & \frac{h}{3}(h)^{3} \sum_{j_{1} ,j_{2}} (|D_{+1}D_{+2}u_{j_{1} +k_{1},j_{2}}-D_{+1}D_{+2}u_{j_{1} ,j_{2}}|^{2} + \\ & & \mbox{}
     |D_{+1}D_{+2}u_{j_{1} ,j_{2}+k_{2}}-D_{+1}D_{+2}u_{j_{1} ,j_{2}}|^{2} \\
& \leq & c_{1}h^{4}(k_{1}(k_{1}-1)\frac{C}{h^{2}} + k_{2}(k_{2}-1)\frac{C}{h^{2}} \\
& \leq & \tilde{c} |\delta|^{2}.
\end{eqnarray*}

\end{proof}

\subsection{Constructing Discrete Approximations}
\label{constructapprox}
We now construct semi-discrete approximations for  LLG.  The method that we adopt is that of finite-differencing on a uniform spatial grid with gridsize $h$.  This method reduces the problem to a system of o.d.e's with the unknowns being $\{\uh{j}\}$, the value of the function at every grid point. We seek a construction that yields a stable scheme.   A natural way to achieve this is to start by discretizing the energy  rather than by discretizing the equation. 
Explicitly, the Dirichlet energy 
can be discretized to yield  the \emph{discrete energy} $E^{h}[u]$:
\[
E^{h}[u^{h}] = \frac{1}{2}h^{2}\sum_{j}\frac{1}{2}(|D_{+\alpha}\uh{j}|^{2} + |D_{-\alpha}\uh{j}|^{2}),
\]
where $j$ indexes over the grid pairs $(j_{1}, j_{2})$.  The discrete analog of $\frac{\delta E[u]}{\delta u}$ is
\[
\frac{\delta E^{h}[u^h]}{\delta \uh{j}} = -\frac{\uh{j+1} - 2\uh{j} + \uh{j-1}}{h^2} = -\Delta^h \uh{j} 
\]
The discretization of the system associated to this energy is  then given by
\begin{equation}
\label{dllg}
\left\{\begin{array}{lcl}
\dt{\uh{j}(t)}   = \nu_{j} \wedge \Delta^{h}\uh{j} +  \alpha(\Delta^{h} \uh{j} + \lambdah{j}\nu_{j})  \\
\uh{j}(t = 0)  =   f_{j}^{h}
\end{array}\right.
\end{equation}
where  $\nu_{j}$ denotes the unit normal vector to ${\mathcal{N}}$ at $u_{j}^{h}(t)$ and the Lagrange's multiplier $\lambdah{j}(t)$ is given by $\lambdah{j}(t) =  -\Laph \uh{j}(t)\cdot \nu_{j}(t)$. 

We also have the discrete equivalent of (\ref{llgform2}):
\begin{equation}
\label{dgeneral2}
\alpha \dt{\uh{j}(t)} - \nu_{j} \wedge \dt{\uh{j}}   =   (1+\alpha^2)(\Delta^{h} \uh{j} + \lambdah{j}\nu_{j}) .
\end{equation}
This discretization yields energy bounds.     In particular, since $\nu_{j}$ is
orthogonal to $\dt{\uh{j}}$, we can multiply $(\ref{dgeneral2})$ by  $\displaystyle \dt{\uh{j}}$,
and sum across indices $j$ to get 
\begin{equation}
\label{freebounds}
\dt{E^{h}[u^{h}(t)]}
       =  (\dt{u^{h}}, \frac{\delta  E[u^{h}]}{\delta u^{h}})_{L_{h}^{2}}\\
       =  - \alpha\|\dt{u^{h}}\|_{L_{h}^{2}}^{2} .
\end{equation}
Since $\alpha > 0$,  we can integrate on any time interval $[0,T]$  to get the discrete energy bound
\[
E^h[u^h[T]] \leq E^h[f^h].
\]

\subsection{Constructing Discrete Fundamental Solutions}
Our method of proof is based on linear estimates.  Although not immediately apparent, the relevant linear equation to consider is the linear damped Sch\"odinger equation.  Since we will also be looking at the harmonic map heat flow, we will also consider just the linear heat equation.  We now construct fundamental solutions for the discrete analogs of these equations.  The fundamental solution of the discrete heat equation is given by  $\Phi^{h}(t) = (\Phi_{j}^{h})_{j}$, where
$\Phi_{j}^{h}$ solves the system
\[
\left\{\begin{array}{lcl}
\dt{\Phi_{j}^{h}(t)} = \Delta^{h}\Phi_{j}^{h} \\
\Phi_{j}^{h}(t = 0)  =   \frac{1}{h^{d}}\delta^{h}
\end{array}\right.
\]
The notation is consistent with the definition of the weighted $h$ norms that is being used and yields 
$\|\frac{1}{h^{d}}\delta^{h}\|_{L_{h}^{2}(\R{d})} = 1$.
An explicit description of the fundamental solution follows from  the observation that the 
discrete Laplacian  $\Delta^{h}u_{j}^{h} = \Delta^{h}u_{j_{1}, j_{2}, \ldots, j_{d}}^{h}$ can be written
as the sum of operators,
\[
 \Delta^{h}u_{j}^{h} = \sum_{i = 1}^{d} 
\frac{1}{h}(D_{+i}- D_{-i})u_{j}^{h}.
\]
As in the continuous case,  the description
of the fundamental solution in higher dimension  rests on the one dimensional problem of solving
\begin{equation}
\label{dkernel}
\left\{\begin{array}{lcl}
\dt{w_{j_{i}}^{h}(t)} = \frac{1}{h}D^{+}w_{j_{i}}^{h} - \frac{1}{h}D^{-}w_{j_{i}}^{h} \\
w_{j_{i}}^{h}(t = 0)  =   \frac{1}{h}\delta_{i}^{h}
\end{array}\right.
\end{equation}
The kernel associated to the operator $\ddt - \frac{1}{h}D^{+}$ is
given by $K_{+}^{h}$, where  
\[
(K_{+}^{h})_{j_{i}}(t)= \left\{ \begin{array}{ll} 
  \frac{1}{h}\frac{(\frac{t}{h^{2}})^{-j_{i}}}{(-j_{i})!} e^{\frac{-t}{h^{2}}} & \mbox{if $j_{i} \leq 0$} \\
 0 & \mbox{otherwise}
\end{array}
\right.,
\]
and that associated with $\frac{\partial}{\partial t} + \frac{1}{h}D^{-}$ is given by $K_{-}^{h}$ where
\[
(K_{-}^{h})_{j_{i}}(t)= \left\{ \begin{array}{ll} 
  \frac{1}{h}\frac{(\frac{t}{h^{2}})^{j_{i}}}{(j_{i})!} e^{\frac{-t}{h^{2}}} & \mbox{if $j_{i} \geq 0$} \\
 0 & \mbox{otherwise}
\end{array}
\right., 
\]
we exploit once again the commutativity of $D^{+}$ and $D^{-}$ to conclude that the
solution to $(\ref{dkernel})$ is given by the discrete convolution 
\begin{equation}
\label{1dkernel}
w^{h}_{j_{i}}(t)  = 
\sum_{k \geq j_{i}}  \frac{1}{h}\frac{(\frac{t}{h^{2}})^{k-j_{i}}}{(k-j_{i})!} e^{\frac{-t}{h^{2}}} 
 \frac{(\frac{t}{h^{2}})^{k}}{(k)!} e^{\frac{-t}{h^{2}}} . 
\end{equation}
The fundamental solution of the discrete heat operator in  $\R{d}$, $\Phi_{j}^{h}$, is then just
the product of the $w_{j_{i}}^{h}$'s over $1 \leq i \leq d$.

Similarly, the fundamental solution of the discrete damped Schr\"odinger equation, denoted by $\U^{\alpha,h} = (\U_{j}^{\alpha,h})_{j}$,  solves the system
\[
\left\{\begin{array}{lcl}
\dt{\U_{j}^{\alpha,h}(t)} = (\alpha + i)\Delta^{h}\U_{j}^{\alpha,h} \\
\U_{j}^{\alpha,h}(t = 0)  =   \frac{1}{h^{m}}\delta^{h}
\end{array}\right.
\]
The construction proceeds exactly as in the case of the heat equation.
The following proposition summarizes these constructions:
\begin{prop}
\label{dkernels}
\hfill
\begin{enumerate}
\item $\Phi^{h} = (\Phi_{j}^{h})_{j}$, the fundamental solution for the $d$-dimensional discrete heat equation, is given by
\begin{equation}
\label{dheatkernel}
\Phi_{j}^{h}(t) = \Pi_{i = 1}^{d}  \sum_{k \geq j_{i}}  \frac{1}{h}\frac{(\frac{t}{h^{2}})^{k-j_{i}}}{(k-j_{i})!} 
 \frac{(\frac{t}{h^{2}})^{k}}{(k)!} e^{\frac{-2t}{h^{2}}} . 
\end{equation}
The solution operator to the linear discrete heat equation for inital value  $f^{h}$ is
given by the discrete convolution   $\Phi^{h}(t) \ast f^{h}$ and  
will be denoted by $\Phi^{h}(t)f^{h}$.
\item $\U^{\alpha,h}= (\U_{j}^{\alpha,h})_{j}$, the fundamental solution for the $d$-dimensional discrete damped Schr\"odinger equation, is given by 
\begin{equation}
\label{dschrokernel}
\U_{j}^{\alpha,h}(t) = \Pi_{i = 1}^{d}  \sum_{k \geq j}  \frac{1}{h}\frac{((\alpha+i)\frac{t}{h^{2}})^{k-j}}{(k-j)!}
 \frac{((\alpha+i)\frac{t}{h^{2}})^{k}}{(k)!} e^{-2(\alpha+i)\frac{t}{h^{2}}} . 
\end{equation}
The solution operator to the linear discrete damped Schr\"odinger equation for inital value  $f^{h}$ is denoted by $\U^{\alpha,h}(t)f^{h}$.   Additionally, denote $\U^{1,h} = \U^{h}$.  
\end{enumerate}
\end{prop}

\subsection{Discrete Linear Estimates}
Having constructed the discrete fundamental solutions, it is possible to derive discrete analogs of
the well known $L^{p}-L^{q}$ estimates.  
\begin{prop} ($\Lh{p}-\Lh{q}$ Estimates)
Let $\K^{h} = \Phi^{h}$ or $\U_{\alpha, h}$, $\alpha > 0$  from Proposition $\ref{dkernels}$.  Then the following
estimates hold:
\mbox{}
\begin{itemize}
\item[(1)] \begin{equation}
\label{de1}
\|\K^{h}(t)f^{h}\|_{\Lh{p}(\R{d})} \leq C\frac{\|f^{h}\|_{\Lh{q}}}{t^{\frac{d}{2}(\frac{1}{q}-\frac{1}{p})}},
\end{equation}
\item[(2)] \begin{equation}
\label{de2}
\|D^{1}\K^{h}(t)f^{h}\|_{\Lh{p}(\R{d})} \leq C\frac{\|f^{h}\|_{\Lh{q}}}{t^{\frac{1}{2} + 
\frac{d}{2}(\frac{1}{q}-\frac{1}{p})}}.
\end{equation}
\end{itemize}
\end{prop}

\begin{proof}
We will  show the derivation of both estimates by taking $\K^{h} = \Phi^{h}$.  Since $\alpha > 0$, the
same derivation holds in the case that $\K^h = \U^{\alpha, h}$.

These estimates will follow upon obtaining suitable $\Lh{1}$ and $\Lh{\infty}$ estimates on 
$\Phi^{h}$ and $D^{1} \Phi^{h}$.  Specifically, if the bounds 
\[
\|\Phi^{h}(t)\|_{L_{h}^{1}} \leq 1, \ \ 
\|\Phi^{h}(t)\|_{L_{h}^{\infty}} \leq \frac{1}{t^{\frac{d}{2}}}
\]
hold, then by Hausdorff-Young we have 
\[
\|\Phi^{h}(t)f^{h}\|_{\Lh{p}} \leq \|\Phi^{h}(t)\|_{\Lh{r}} \|f^{h}\|_{\Lh{q}}, 
\]
where $1 + \frac{1}{p} = \frac{1}{r} + \frac{1}{q}$.  By interpolation, 
$\|\Phi^{h}(t)\|_{\Lh{r}} \leq (\|\Phi^{h}(t)\|_{\Lh{1}})^{1-\theta}(\|\Phi^{h}(t)\|_{\Lh{\infty}})^{\theta}$
for $\theta = 1-\frac{1}{r} =1-( 1 + \frac{1}{p} - \frac{1}{q}) = \frac{1}{q} - \frac{1}{p}$, so 
$\|\Phi^{h}(t)\|_{\Lh{r}} \leq \frac{1}{t^{\frac{d}{2}(\frac{1}{q} - \frac{1}{p})}}$, from which
 (1) follows.   Similarly,  (2)  follows from the estimates 
\[
\|D^{1}\Phi^{h}(t)\|_{\Lh{1}} \leq \frac{1}{t^{\frac{1}{2}}}, \ \ 
\|D^{1}\Phi^{h}(t)\|_{\Lh{\infty}} \leq \frac{1}{t^{\frac{1}{2} + \frac{d}{2}}}.
\]

\begin{itemize} 
\item  To show  $\|\Phi^{h}(t)\|_{\Lh{1}} \leq 1$,  it suffices to  get an estimate on the  kernel of the one dimensional
problem,
$w^{h}_{j_{i}}(t)$ given by $(\ref{1dkernel})$ which is a discrete convolution of two functions in $L_{h}^{1}$.  Since the $L_{h}^{1}$ bound of 
a convolution of two $\Lh{1}$ functions can be bounded by the product of the two $\Lh{1}$ bounds, we have
\begin{eqnarray*}
\|w^{h}_{j_{i}}(t)\|_{\Lh{1}} &  =  &\|\sum_{k \geq j_{i}}  
\frac{1}{h}\frac{(\frac{t}{h^{2}})^{k-j_{i}}}{(k-j_{i})!} e^{\frac{-t}{h^{2}}} 
 \frac{(\frac{t}{h^{2}})^{k}}{(k)!} e^{\frac{-t}{h^{2}}} \|_{\Lh{1}} \\
& \leq & \|e^{-\frac{2t}{h^{2}}} \sum_{k \geq j_{i}}  
\frac{1}{h}\frac{(\frac{t}{h^{2}})^{k-j_{i}}}{(k-j_{i})!}  \frac{(\frac{t}{h^{2}})^{k}}{(k)!}\|_{\Lh{1}} \\
& \leq & e^{-\frac{2t}{h^{2}}}\|\frac{1}{h} \frac{(\frac{t}{h^{2}})^{k}}{(k)!}\|_{\Lh{1}}^{2} 
\end{eqnarray*}
\begin{eqnarray*}
& \leq & e^{-\frac{2t}{h^{2}}} e^{\frac{2t}{h^{2}}} \\
& \leq & 1.
\end{eqnarray*}
Since $w^{h}_{j_{i}}(t)$ is a convolution of two $\Lh{1}$ functions (in $j_{i}$) , the 
$\Lh{\infty}$ bound of $w^{h}_{j_{i}}(t)$ will be the product of the $\Lh{\infty}$ bound of any
one of the convolved pair with the  $\Lh{1}$ bound of the other.  
Using Stirling's formula, and re-indexing with $k$, we have 
\begin{eqnarray*}
\|w^{h}_{k}(t)\|_{\Lh{\infty}} & \leq &
 \|\frac{1}{h}\frac{(\frac{t}{h^{2}})^{k}}{k!} e^{\frac{-t}{h^{2}}}\|_{\Lh{\infty}} \\
 & \leq &  \|\frac{1}{h}\frac{(\frac{t}{h^{2}})^{k}}{C_{k} k^{k}e^{-k}\sqrt{k}} 
       e^{\frac{-t}{h^{2}}}\|_{\Lh{\infty}} \\
 & \leq & C\|\frac{1}{h\sqrt{k}} (\frac{t}{kh^{2}})^{k}e^{ -\frac{t}{h^{2}}+ k}\|_{\Lh{\infty}} \\.
\end{eqnarray*}
Writing $x = \frac{t}{h^{2}}$, and letting  $f(k) = (\frac{x}{k})^{k}e^{ -x+ k}$, we have
\begin{eqnarray*}
\frac{\partial f}{\partial k} & = & ((\frac{x}{k})^{k})'e^{-x+k} + (\frac{x}{k})^{k}e^{-x+k} \\
      & = & (\frac{x}{k})^{k}(\ln{\frac{x}{k}} - 1)e^{-x+k} + (\frac{x}{k})^{k}e^{-x+k}  \\
      & = & e^{-x+k}(\frac{x}{k})^{k}\ln{\frac{x}{k}}
\end{eqnarray*}
Thus $f$ is maximum when $\frac{\partial f}{\partial k}$ equals $0$ which occurs when $k = x = \frac{t}{h^{2}} \Rightarrow \frac{1}{h} = \sqrt{\frac{k}{t}}$.
Returning to the $\Lh{\infty}$ bound of $w^{h}_{k}(t)$, we have
\begin{equation}
\label{winfty}
\|w^{h}_{k}(t)\|_{\Lh{\infty}}  \leq  C\|\frac{1}{h\sqrt{k}} f(k)\|_{\Lh{\infty}} 
  \leq  C\|\frac{1}{h\sqrt{k}}\|_{\Lh{\infty}}f(x) =  C\frac{1}{\sqrt{t}}
\end{equation}
as needed.  It follows that $\|\Phi^{h}(t)\|_{\Lh{\infty}} \leq \frac{C}{t^{\frac{d}{2}}}$.

\item To get gradient bounds on the fundamental solution, we express 
\[
D_{+\alpha}\Phi^{h}(t) =  D_{+\alpha} w_{j_{\alpha}}^{h}(t)  \Pi_{i \neq \alpha} w_{j_{i}}^{h}(t)
\]
For $\Lh{\infty}$ bounds on the gradient, we use the fact that the gradient will
only 'hit' one term in the convolved pair and that the other is in $\Lh{1}$.  Again, let $x = \frac{t}{h^{2}}$, and $f(k) = (\frac{x}{k})^{k}e^{ -x+ k}$.  Since $f(k)$ achieves 
its maximum at $k = x$,
\[
|D_{+\alpha}w_{k}^{h}(t)| \leq |\frac{1}{h^{2}} f(k) (\frac{t}{h^{2}}\frac{1}{j_\alpha+1}-1)| \leq f(x)\frac{k}{k+1} \frac{1}{t} \leq C\frac{1}{t}.
\] 
Using the $\Lh{\infty}$ bound given by ($\ref{winfty}$) for each of the other  $d-1$ $w_{j_{i}}^{h}(t)$ terms in the product,  we get
\[
\|D_{+\alpha}\Phi^{h}(t)\|_{L_{h}^{\infty}} \leq \frac{C}{t^{1 + \frac{d-1}{2}}} = \frac{C}{t^{\frac{1}{2} + \frac{d}{2}}}.
\]
Therefore, 
\[
\|D^1 \Phi^{h}(t)\|_{L_{h}^{\infty}} \leq  \frac{C}{t^{\frac{1}{2} + \frac{d}{2}}}.
\]
Employing again the fact that the $\Lh{1}$ bound of a convolution can be bound by the product of each
of the $\Lh{1}$ bounds  in the convolved pair (and hence 
$\|D_{+\alpha}w_{j_{\alpha}}^{h}(t)\|_{L_{h}^{1}} \leq \frac{C}{h} \leq \frac{C}{t^{\frac{1}{2}}}$ for each $\alpha$),  and that each of the $w_{j_{i}}^{h}(t)$ is an $\Lh{1}$ function ,  we conclude that 
\[
\|D^1 \Phi^{h}(t)\|_{L_{h}^{1}} \leq  
C \sup_{1 \leq i \leq d} \|(D_{+i}w_{j_{i}}^{h}(t)) \Pi_{l \neq i} w_{j_{l}}^{h}(t)\|_{\Lh{1}} 
                    \leq  \frac{C}{t^{\frac{1}{2}}}.
\]
\end{itemize}
\end{proof}

\section{Harmonic Map Heat Flow}
\label{hhp}
The first problem on which we apply the machinery that has been developed so far is the Cauchy problem for the  harmonic map heat flow problem in two dimensions.  This is the problem of finding $u$, given $f : \R{2}\rightarrow \mathcal{N}$ such that
\begin{equation}
\label{ceq}
\left\{\begin{array}{lcl}
\dt{u} & = &  \Delta u + (\nabla \nu(u) \cdot \nabla u) \nu(u) \\
u(0) & = &  f
\end{array}\right..
\end{equation}
$\mathcal{N}$, as usual, is taken to be a compact smooth hypersurface.
The semi-discrete approximation for this equations follows the construction in Section \ref{constructapprox} and is given by
\begin{equation}
\label{deq}
\left\{\begin{array}{lcl}
\dt{\uh{j}(t)} = \Delta^{h}\uh{j}(t) + \lambdah{j}(t)\nu_{j}(t) \\
u_{j}^{h}(0)  =   f_{j}^{h}
\end{array}\right.,
\end{equation}
where  $\nu_{j} = \nu(\uh{j})$  and  the Lagrange's multiplier $\lambdah{j}(t)$ is given by 
\[
\lambdah{j}(t) =  -\Laph \uh{j}(t)\cdot \nu_{j}(t). 
\]
In what follows,  the letter $C$ will designate a generic constant that may depend on $\mathcal{N}$ but is independent of the gridsize $h$ or any solution unless specified explicity like $C(E^{h}[f^{h}])$.  

\subsubsection*{Estimates on $\lambdah{j}$}
\label{secLambda}
For the purpose of attaining bounds on $\uh{j}$ solving $(\ref{deq})$, it will be useful to preface the 
technical detail with several observations concerning  $\lambdah{j}$.   In the continuous
setting, there is an explicit expression for the Lagrange's multiplier due to the fact that $\partial_{i}u \in T_u \mathcal{N}$.  
This permits the following manipulation:
\[
-\Delta u \cdot \nu(u) = \partial_{i}(-\partial_{i} u \cdot \nu(u)) + (\partial_{i}u \cdot \partial_{i} \nu(u)) = \nabla \nu(u) \cdot \nabla u. 
\]
Thus, 
\begin{equation}
\label{clambdabound}
|\lambdah{j}| = |(\nabla \nu(u) \cdot \nabla u)| \leq C| \nabla u|^{2}.
\end{equation}
In the discrete setting, difference derivatives lose their tangency which robs $\lambdah{j}$ of as clean an expression.  Writing  $a_{j}^{i} = D_{+k}\uh{j}\cdot \nu_{j}$, we can attempt the same manipulation as in the continuous case to get the following expression for $\lambdah{j}$:
\[
\lambdah{j} = D_{-i}a_{j}^{i} - (D_{-i}u \cdot D_{-i}\nu_{j})
\]
The first observation is that even though this expression  has a local defect, we can
still establish a discrete analog of (\ref{clambdabound}).

\vspace{0.1in}
\noindent \underline{Observation 1}
\begin{equation}
\label{dlambdabound}
|\lambdah{j}| \leq C(| D_{+i}u_{j}|^{2} +| D_{-i}u_{j}|^{2}),
\end{equation}

\noindent  Expressing 
\[
\lambdah{j} = -\frac{1}{h}\{(D_{+i}\uh{h} - D_{-i}\uh{j}) \cdot \nu_{j}\}
\]
we need only show that for fixed index $k$, 
\[
|\frac{1}{h}(D_{+k}\uh{j} \cdot \nu_{j})| = |(\uh{j_{k}+1}- \uh{j_{k}})\cdot  \nu_{j}| \leq C|D_{+k}\uh{j}|^{2}.
\]
Fix some $\delta > 0$.  Then if $|\uh{j_{k+1}}-\uh{j_{k}}| \geq \delta$,] we have 
\[
|(\uh{j_{k+1}}-\uh{j_{k}})\cdot \nu_{j}| \leq \frac{|\uh{j_{k+1}}-\uh{j_{k}}|}{\delta}|\uh{j_{k+1}}-\uh{j_{k}}|
\]
in which case we can take  $C = \frac{1}{\delta}$.  For the case  $|\uh{j_{k+1}}-\uh{j_{k}}| < \delta$,
the desired inequality is  a direct consequence of the compactness
of ${\mathcal{N}}$.  Namely, since  ${\mathcal{N}}$ is  smooth, we can express ${\mathcal{N}}$ locally as the graph
of some quadratic function, so there is a ball about every point such that ${\mathcal{N}}$ behaves quadractically.
Due to the compactness of 
${\mathcal{N}}$, there is a finite subcover of these
balls.  If we choose $\delta$ to be the minimum of these radii, and to specify henceforth that the gridsize $h$ is to
be smaller than $\delta$, then the interpretation of 
$(\uh{j_{k+1}}-\uh{j_{k}})\cdot \nu_{j}$ as the projection of $\uh{j_{k+1}}-\uh{j_{k}}$
onto $\nu_{j}$ immediately yields the given inequality.  

\vspace{0.2in} 
\noindent  \underline{Observation 2}

\noindent For $\delta > 0 $, choose $h$ such that for $k = 1, 2$, 
$|\uh{j_{k}+1}-\uh{j}| < \delta, \mbox{    } |\uh{j} - \uh{j_{k}-1}| < \delta$. 
Then 
\begin{equation}
\label{observation2}
|D^{1}\lambdah{j}| =  O(|D^{1}\uh{l}|^{3}, |D^{1}\uh{l}||D^{2}\uh{l}|),
\end{equation}
where $l$ indexes over next to nearest neighbors of $j$. 
\vspace{0.1in}

\noindent This follows from the calculation
\begin{eqnarray*}
D_{+k}\lambdah{j}
& = & D_{+k}(\Laph \uh{j}\cdot  \nu_{j}) \\
& = & \frac{1}{h}(D_{+i}(D_{+k} \uh{j} \cdot \nu_j) - D_{-i}(D_{+k} \uh{j} \cdot \nu_j))\\
& = & \frac{C}{h}((D_{+i}- D_{-i})(h D_+k \uh{j})^2 ) \\
& = & O(|D^{1}\uh{l}|^{3}, |D^{1}\uh{l}||D^{2}\uh{l}|),
\end{eqnarray*}
$l$ indexes over next to nearest neighbors of $j$. 

\subsection{A Global Smooth Solution for Small Energy Data}
The first result that will be shown is that under the assumption that
the initial data has small total energy, there exists a smooth global solution.  While the result in 
this section can be seen as a corollary to the proof of the finite energy data case,  this assumption permits the `fastest' path to  our goal .  The process is laid out here to fix ideas and to motivate the line of reasoning underlying some of the extra steps needed in the next section.  The following result will be proved:

\begin{theorem}
\label{smallhhp}
There is a constant $\epsilon_{0}$ such that for $f^{h} \in H_h^{1}(\R{2})$ and  $\sup_h E^{h}[f^{h}] < \epsilon_{0} $, the associated sequence of $H^{1}$-interpolants $\{p_{h}\}$ converges strongly in $H^{1}$ to $p$, which solves $(\ref{ceq})$ in the sense of distributions.
\end{theorem}
 
\begin{proof}
The proof of this result is a two step process of attaining appropriate bounds on the sequence $u^h$ and using these bounds to construct interpolants that converge to the solution claimed.
\subsubsection*{Step 1:  Attaining Bounds}
The only quantities on which we have bounds immediately are energy and $\displaystyle \dt{u^{h}}$.  Namely
\begin{equation}
\label{energybound}
E^h[u^{h}(T)] + \|\dt{u^{h}}\|_{L^2([0,T], L_{h}^{2})}^{2}  \leq E^h[f^h].
\end{equation}
Define 
\[
\mathcal{F}^h = \{u^h  = \{\uh{j}\}_j|  t^{1/2}\|D^{2} u^{h}(t)\|_{\Lh{2}} < \infty \mbox{     }\forall t \}.
\]
The higher derivative bounds that we aim for here are uniform-in-$h$ bounds in  (i) $L^2(\dot{H}_h^2)$, (ii) $\mathcal{F}^h$ and (iii) Lip$(\dot{H}_h^1)$.
\begin{itemize}
\item [(i)] \fbox{$L^2(H_h^2)$}
To get bounds on terms involved in $\|D^{2}u^{h}\|_{\Lh{2}}$, it suffices to bound $\|\Delta^{h} u^{h}\|_{\Lh{2}}$ since the discrete Laplacian controls all second difference derivatives.  
Using Observation 1 and (\ref{sob-int}) for $j = 1, r = 4, p = q = 2$,  and the global energy bound (\ref{energybound}), 
\begin{eqnarray*}
\|\Delta^h u^h\|_{\Lh{2}} 
& \leq &  C(\|\partial_t u^h(t)\|_{\Lh{2}} + (h^2\Sigma_j |\lambda_j^h|^2)^{\frac{1}{2}}) \\
& \leq & C(\|\partial_t u^h(t) \|_{\Lh{2}} + \||D^1 u^h(t)|^2\|_{\Lh{2}}) \\
& \leq & C(\|\partial_t u^h(t) \|_{\Lh{2}} + \|D^1 u^h(t)\|_{\Lh{4}}^2) \\
& \leq & C(\|\partial_t u^h(t) \|_{\Lh{2}} + \|D^1 u^h(t)\|_{\Lh{2}}\|D^2 u^h(t)\|_{\Lh{2}})\\
& \leq & C(\|\partial_t u^h(t) \|_{\Lh{2}} + E^h[f^h]^{\frac{1}{2}}\|D^2 u^h(t)\|_{\Lh{2}})\\
& \leq & C(\|\partial_t u^h(t) \|_{\Lh{2}} + \epsilon_0^{\frac{1}{2}}\|D^2 u^h(t)\|_{\Lh{2}})
\end{eqnarray*}
Now choosing $\epsilon_0$ sufficiently small, we can absorb the $C(\epsilon_0^{\frac{1}{2}})\|D^2 u^h\|_{\Lh{2}}$
term into the left hand side to get
\[
\|\Delta^h u^h\|_{\Lh{2}} \leq C(\|\partial_t u^h(t) \|_{\Lh{2}} 
\]
Now integrating in time, we can make use of ($\ref{energybound}$) again to conclude that 
\[
\|\Delta^h u^h\|_{L^2(\Lh{2})} \leq C(\|\partial_t u^h(t) \|_{L^2(\Lh{2})} \leq  C(E^h[f^h]^{1/2}) \leq C(\epsilon_0).
\]

\item[(ii)] \fbox{$\mathcal{F}^h$} For each index $j$, $u_j^h$ satisfies
\[
 \uh{j}(t)  =  \Phi_{j}^{h}(t)f_{j}^{h} + \int_{0}^{t}\Phi_{j}^{h}(t-s)|\lambda_{j}^h(s)\nu(\uh{j}(s))|  ds.
\] 
To get  bounds involving $D^{2} u^{h}$ bounds, we can take difference derivatives of the equation
and and then gain a derivative by using the estimates on the derivative of the fundamental solution.
Using the estimate given by $(\ref{de2})$, once with $p = q = 2$ and another with $p = 2, q = p' = 4/3$, 
\begin{eqnarray*}
\|D^{2}u^{h}(t)\|_{\Lh{2}}
   & \leq & \frac{\|D^1 f^{h}\|_{\Lh{2}}}{t^{1/2}} + 
            C\int_{0}^{t} \frac{1}{(t-s)^{3/4}}\|D^{1}(\lambda^{h}(s)\nu(u^h(s)))\|_{\Lh{\frac{4}{3}}} ds.
\end{eqnarray*}
Observation 1 and 2 on $\lambdah{j}$ together with the smoothness of $\nu$ yield the inequality
\[
D^{1}(\lambdah{j} \nu_{j})
 \leq C|(D^{1}\lambdah{j}| + |\lambdah{j}||D^{1}\nu_{j}|)
 \leq  C(|D^{1}u_{l}|^{3} +  |D^{1}u_{l}||D^{2}u_{l}|), 
\]
where $l$ indexes over next to nearest neighbors of $j$. 
We can now apply  ($\ref{sob-int}$) in addition
to repeated use  of the discrete Holder and Minkowski inequalities to get
\[
\|D^{1}(\lambda^{h}\nu)\|_{\Lh{\frac{4}{3}}}
 \leq   C(\|D^{1}u^{h}\|_{L_{h}^{4}}\|D^{2}u^{h}\|_{\Lh{2}} + \|D^{1}u^{h}\|_{\Lh{4}}^{3}) . 
\]
Proceeding,    
\begin{eqnarray*}
\lefteqn{t^{1/2}\|D^{2}u^{h}(t)\|_{\Lh{2}}} \\
   & \leq & C(\|D^{1} f^{h}\|_{\Lh{2}}) + Ct^{1/2}\int_{0}^{t} \frac{1}{(t-s)^{3/4}}\{\|D^{1}u^{h}(s)\|_{\Lh{4}}
            \|D^{2}u^{h}(s)\|_{\Lh{2}}  + \|D^{1}u^{h}\|_{\Lh{4}}^{3}\} ds \\ 
   & \leq &  CE^h[f^h]^{1/2}+ 
            C t^{1/2}\int_{0}^{t} \frac{1}{(t-s)^{3/4}}
           \|D^{1}u^{h}(s)\|_{\Lh{2}}^{1/2}\|D^{2}u^{h}(s)\|_{\Lh{2}}^{3/2} + \\
           &  & \mbox{} + \|D^{1}u^{h}(s)\|_{\Lh{2}}^{3/2}\|D^{2}u^{h}(s)\|_{\Lh{2}}^{3/2}ds\\
   & \leq &   CE^h[f^h]^{1/2}+ 
            C(\sup_{t \geq 0} E^{h}[u^{h}(t)]) t^{1/2}\int_{0}^{t} \frac{1}{(t-s)^{3/4}}
             \|D^{2}u^{h}(s)\|_{\Lh{2}}^{3/2} ds\\   
   & \leq &  C\epsilon_0^{1/2} + 
            C(\epsilon_{0}) t^{1/2}\int_{0}^{t} \frac{1}{(t-s)^{3/4}}
           \frac{1}{s^{3/4}} (s^{1/2} \|D^{2}u^{h}(s)\|_{\Lh{2}})^{3/2} ds.
\end{eqnarray*}  
Letting $y(t) = \sup_{0 \leq s \leq  t}  s^{1/2}\|D^{2}u^{h}(s)\|_{\Lh{2}}$,
\[
 y(t) \leq  \epsilon_{0}^{1/2} +  c(\epsilon_{0}) y(t)^{3/2} t^{1/2} \int_{0}^{t} \frac{1}{(t-s)^{3/4}}
           \frac{1}{s^{3/4}} ds.
\]
Noting  that 
\[
t^{1/2}\int_{0}^{\frac{t}{2}}\frac{1}{(t-s)^{3/4}} \frac{1}{s^{3/4}} ds 
   \leq t^{1/2}(\frac{2}{t})^{3/4}\int_{0}^{\frac{t}{2}}\frac{1}{s^{3/4}} ds 
   \leq C
\]
and 
\[
t^{1/2} \int_{\frac{t}{2}}^{t}\frac{1}{(t-s)^{3/4}} \frac{1}{s^{3/4}} ds 
  \leq t^{1/2} (\frac{2}{t})^{3/4}\int_{\frac{t}{2}}^{t}\frac{1}{(t-s)^{3/4}} ds
\leq C
\]
we have that $y(t)$ satisfies
\[
y(t) \leq C \epsilon_{0}^{1/2} +  C(\epsilon_{0}) y(t)^{3/2}
\] 
for all $t \in [0, \infty)$.  The function $L(y) = y- C(\epsilon_0) y^{3/2}$ has one root at $y = 0$ and another positive root, reaching its maximum at $y = \frac{4}{9C(\epsilon_0)^2}$.  Now we can choose 
$\epsilon_{0}$ so that $C\epsilon_0^{1/2}  < L(\frac{4}{9C(\epsilon_0)^2})$.  Since $t \mapsto y(t)$ is continuous and  $y(0) = 0$, we can conclude that $y(t)$ is bounded, uniform in $h$.

\item[(iii)] \fbox{Lip$(\dot{H}_h^1)$} The energy bound (\ref{energybound}) immediately implies that $u^h \in L^{\infty}(\dot{H}_h^1)$.  We can do better using the
bounds that we just attained in $\mathcal{F}^h$.  For $t_1 < t_2 $, let $t^* = t_1 + \frac{t_2 -t_1}{2}$.   The bound we seek is a consequence of the following:
\begin{eqnarray*}
\lefteqn{\|D^1 (u^h(t_2) - u^h(t_1))\|_{\Lh{2}}}\\
& \leq & C(\int_{t_1}^{t^*} \frac{1}{|t_1 - s|^{1/2}}\|D^1 u^h(s)|_{\Lh{4}}^2 ds+ \int_{t^*}^{t_2} \frac{1}{|t_2 - s|^{1/2}}\|D^1 u^h(s)|_{\Lh{4}}^2 ds)\\
& \leq & C(\int_{t_1}^{t^*} \frac{1}{|t_1 - s|^{1/2}}\|D^1 u^h(s)|_{\Lh{2}} \|D^2 u^h(s)|_{\Lh{2}} ds+  \\
& & \mbox{} +  \int_{t^*}^{t_2} \frac{1}{|t_2 - s|^{1/2}}\|D^1 u^h(s)|_{\Lh{2}} \|D^2 u^h(s)|_{\Lh{2}}ds)\\
& \leq & C(\epsilon_0)(\int_{t_1}^{t^*} \frac{1}{|t_1 - s|^{1/2}}\frac{1}{s^{1/2}} ds+  \int_{t^*}^{t_2} \frac{1}{|t_2 - s|^{1/2}}\frac{1}{s^{1/2}} ds)\\
& \leq & C(\epsilon_0)((t^*-t_1)+ (t_2-t^*) )\\
& = & C(\epsilon_0)(t_2-t_1)
\end{eqnarray*}
\end{itemize}

\subsubsection*{Step 2: Convergence of the Interpolants}
\noindent We have shown that there exists, for each $h$, sequences of solutions $u^h = \{\uh{j}(t)\}$ such
that 
\[
u^h \in L^2(\dot{H}_h^2) \cap \mathcal{F}^h \cap \rm{Lip}(\dot{H}_h^1)
\]
with uniform bounds in $h$:
\[
\sup_h( \|u^h\|_{L^2(\dot{H}_h^2)} +\sup_{t}  t^{1/2}\|D^{2} u^{h}(t)\|_{\Lh{2}} + \sup_{t_1 \neq t_2} \frac{\|u^h(t_2)-u^h(t_1)\|_{\dot{H}_h^1}}{|t_2 -t_1|}) < C(\epsilon_0)
\] 
Associate with each sequence $\{u_{j}^{h}(t)\}$ the $H^{1}$ interpolant $p_{h}(t)$.  Then for $\phi \in C_{0}^{\infty}(\Rplus \times \R{2})$, 
we will show the following statement:
\begin{equation}
\label{phintegral}
-\int_0^{\infty} \int_{\R{2}} p_h \cdot \partial_t \phi  + p_h \cdot \Delta \phi + (\nabla p_h \cdot \nabla \nu(p_h))\nu(p_h) \cdot \phi  = \int_{\R{2}}  p_h(0)\cdot \phi(0)  + O(h)
\end{equation}
Recall that, for $(x,y) \in \Box_{j}$, $j = (j_{1}, j_{2})$, the $H^{1}$ interpolant is given by 
\[
p_{h}(x,y,t) = \uh{j}(t) + D_{+1}\uh{j}(t)(x-j_{1} h) + D_{+2}\uh{j}(t)(y-j_{2}h) + D_{+1}D_{+2}\uh{j}(t)(x-j_{1} h)
(y-j_{2}h)).
\]
For $j = (j_1, j_2)$, the notation $\int_{\Box_j}$ will denote $\int_{j_{1}h}^{(j_{1}+1)h}\int_{j_{2}h}^{(j_{2}+1)h}$.  We will also slightly abuse the notation $u^h$ (which is only defined on grid points) to
play the dual role of the function which takes the value of $u_j^h$ in $\Box_j$.
Since $u^h \in L^2(\dot{H}_h^2) \cap L^{\infty}(\dot{H}_h^1)$, the following are immediate:
\begin{itemize}
\item[(i)] $\displaystyle \int_0^{\infty} (\int_{\R{2}} p_h \cdot \partial_t \phi - \sum_j \int_{\Box_j} \uh{j} \cdot \partial_t \phi )= O(h)$
\item[(ii)] $\displaystyle \int_0^{\infty} ( \int_{\R{2}} p_h \cdot \Delta \phi -  \sum_j \int_{\Box_j} \uh{j} \cdot \Delta \phi  )= O(h)$
\end{itemize}
From (i),
\[
 \int_0^{\infty} \sum_j \int_{\Box_j} \uh{j} \cdot \partial_t \phi = - \int_0^{\infty} \sum_j \int_{\Box_j} \partial_t \uh{j} \cdot \phi  +  \sum_j \int_{\Box_j} f^h(0) \cdot \phi(0)
\]
Using the equation (\ref{deq}) for $\partial_t \uh{j}$, what remains to be shown is that 
\begin{itemize}
\item[(iii)]  $\displaystyle \int_0^{\infty} \sum_j \int_{\Box_j}  \Delta^h \uh{j} \cdot \phi - p_h \cdot \Delta \phi = O(h)$
\item[(iv)] $\displaystyle \int_0^{\infty} \sum_j \int_{\Box_j}  \lambda_j^h \nu_j \cdot \phi + (\nabla p_h \cdot \nabla \nu(p_h))\nu(p_h) \cdot \phi  = O(h)$
\end{itemize}
To show (iii), it suffices to show that 
\[
\int_{0}^{\infty} \sum_j \int_{\Box_j} \Delta^h \uh{j} \cdot \phi = \int_{0}^{\infty} \sum_j \int_{\Box_j} \uh{j} \cdot \Delta \phi + O(h)
\] 
\begin{eqnarray*}
\int_{0}^{\infty} \sum_j \int_{\Box_j} \Delta^{h}\uh{j} \cdot \phi 
& = & \int_{0}^{\infty} \sum_j \int_{\Box_j} \frac{D_{+i}\uh{j}-D_{-i}\uh{j}{h}}{h} \cdot \phi\\
& = & - \int_{0}^{\infty} \sum_j \int_{\Box_j} D_{+i}\uh{j}\cdot \frac{\phi(x)-\phi(x_i+h)}{h} \\
& = &  \int_{0}^{\infty} \sum_j \int_{\Box_j} \uh{j}\cdot \frac{\phi(x_i +h) - 2\phi(x_i) + \phi(x_i -h)}{h^2}\\
& = &  \int_{0}^{\infty} \sum_j \int_{\Box_j} \uh{j}\cdot \Delta \phi + O(h)
\end{eqnarray*}
What remains is  (iv).  Associate with  $\phi$  the step function $\phis = \sum_{j} \phi_{j}\chi_{\Box_{j}}$.   We treat this term by splitting the integral into two parts:
\[
I = \int_{0}^{\infty} \sum_j \int_{\Box_j} (\lambdah{j} \nu_{j} + (\nabla p_{h}\cdot \nabla\nu(p_{h}))\nu(p_{h}))\cdot (\phi-\phis)
\]
and
\[
II = \int_{0}^{\infty} \sum_j \int_{\Box_j} (\lambdah{j} \nu_{j} + (\nabla p_{h}\cdot \nabla\nu(p_{h}))\nu(p_{h})) \cdot \phis
\]
Since  $u^{h} \in L^{2}(\dot{H}_h^2)$, $ I = O(h) $.  The reason behind the splitting is really to negotiate the local defect regarding 
$\lambdah{j}$.  Fixing $j$, $\lambdah{j}$ does not come close to cancelling out the relevant terms
in $\nabla p_{h}(jh)\cdot \nabla\nu(p_{h}(jh))$.  However, just as it was seen in Observation 1 that the local defect did not prevent us from attaining
a global statement that was satisfactorily close to the continuous case, we can also show that by considering the
sum of the two terms in $II$ as two sums, that the terms that result have a satisfactory cancellation.
The only thing that was preventing us before from switching  the spatial sum and integral was $\phi$.  Replacing $\phi$ by $\phis$ resolves this immediate problem.   
Keeping in mind that the terms in $\nabla p_{h}(jh)\cdot \nabla\nu(p_{h}(jh))$ that require cancellation are 
of the form $D_{+i}\uh{j}\cdot (\nu_{i})_j$, 
we consider the sum $\sum_{j} \lambdah{j}\nu_{j} \cdot \phih_{j}$.   Letting $a_j^i = (D_{+k} \uh{j}\cdot \nu_j)$,

\begin{eqnarray*}
\sum_{j} \lambdah{j}\nu_{j} \cdot \phih_{j}
& = & \sum_{j} D_{-i}a_{j}^{i} - (D_{-i}\uh{j} \cdot D_{-i}\nu_{j}) \nu_{j}\cdot  \phih_{j}\\
& = &  \sum_{j} D_{-i}a_{j}^{i} \cdot \phih_{j} - \sum_{j} (D_{+i}\uh{j} \cdot D_{+i}\nu_{j}) \nu_{j_{i}+1} \cdot \phih_{j_{i}+1}\\
& = &  - \sum_{j} (D_{+i}\uh{j} \cdot D_{+i}\nu_{j}) \nu_{j_{i}+1} \cdot  \phih_{j_{i}+1}
\end{eqnarray*}
This term is already in a form that the cancellation is apparent;  the local defect in $\lambdah{j}$ has essentially
been transferred to a local defect in $\nu$ and  $\phi$, but these we can handle because they are smooth.
\begin{eqnarray*}
II &=& \int_{0}^{\infty}\int_{\Box_j}  \sum_{j} \phi_{j} \cdot \{ \lambdah{j} \nu_{j}  + D_{+i}\uh{j})(\nu_{x})_{j}\nu(p_{h})\} + O(h)\\
& = &  \int_{0}^{\infty} h^{2}  \sum_{j} \lambdah{j}\nu_{j}\cdot \phi_{j}+ ( D_{+i}\uh{j}\cdot (\nu_{i})_{j})\nu_{j}\cdot \phi_{j}  + O(h)\\
& = &  \int_{0}^{\infty} h^{2}\{ - \sum_{j} (D_{+i}\uh{j} \cdot D_{+i}\nu_{j}) \nu_{j_{i}+1} \cdot  \phi_{j_{i}+1} +   \sum_{j} (D_{+i}\uh{j} \cdot (\nu_{i})_{j})\nu_{j}\cdot \phi_{j} + O(h)\\
& = & O(h).
\end{eqnarray*}
(i) - (iv) together entail (\ref{phintegral}).
With this expression, we can split the time interval into subintervals $[0, \delta)$ and $[\delta, \infty)$.  On $[\delta,  \infty) $, we have uniform $ \dot{H}_h^2$ which implies from Proposition ${\ref{phlemma}}$, that there exists a subsequence $\{h_{k}\}$ such that $p_{h_{k}}(t) \stackrel{H^{1}}{\rightarrow} p(t) $. 
For any $h_k$, then, we can use in addtion (i) - (iv)  and the fact that $u^h \in \rm{Lip}(\dot{H}_h^1)$ to show that 
\begin{eqnarray*}
\lefteqn{ \int_{0}^{\infty}\int_{\R{2}}p_{h_k} \cdot \partial_t \phi  + p_{h_k} \cdot \Delta \phi + (\nabla p_{h_k} \cdot \nabla \nu(p_{h_k}))\nu(p_{h_k}) \cdot \phi  }\\ 
& = & (\int_{0}^{h_k} + \int_{h_k}^{\infty} )\int_{\R{2}} p_{h_k} \cdot \partial_t \phi  + p_{h_k} \cdot \Delta \phi + (\nabla p_{h_k} \cdot \nabla \nu(p_{h_k}))\nu(p_{h_k}) \cdot \phi \\
& = & C(\epsilon_0)h_k  + O(h_k) + \int_{h_k}^{\infty} \int_{\R{2}} p_{h_k} \cdot \partial_t \phi  + p_{h_k} \cdot \Delta \phi + (\nabla p_{h_k} \cdot \nabla \nu(p_{h_k}))\nu(p_{h_k}) \cdot \phi \\
& \rightarrow &  \int_{0}^{\infty}\int_{\R{2}} p \cdot \partial_t \phi  + p \cdot \Delta \phi + (\nabla p \cdot \nabla \nu(p))\nu(p) \cdot \phi
\end{eqnarray*}
from which we can conclude that $p$ solves (\ref{ceq}) in the sense of distributions.
\end{proof}

\subsection{A Partially Regular Solution for Finite Energy Data}
Now we consider the discrete system  only assuming  $H^{1}$ data.  While we do not have total small energy, our initial data is
still of finite energy.  Locally, then, the energy of the initial data can be made small.  What is needed is a localized version of the
``small energy $\Rightarrow$ regularity'' statement that was gotten in the last section.  Achieving such
a statement will allow us to estimate the size of the concentration sets which in turn
will allow us to extend the associated interpolants across these sets in such a way that they will converge as 
$h \rightarrow 0$ in the sense of distributions to a solution of the continuous problem.   We proceed  now to show the following result:
\begin{theorem}
\label{fehhpthm1}
For any $u_{0} \in H^{1}(\R{2}, \mathcal{M})$, there exists a global, partially regular solution to (\ref{ceq}), smooth away from a locally finite $2$-dimensonal  Hausdorff measure set with respect to the parabolic metric.
\end{theorem}

\noindent The core statement to be shown is:

\begin{lemma}
\label{fehhpthm2}
There is an $\epsilon_{0} > 0$ such that if $E^{h}[u^{h}(t_{0}); B_{\overline{R}}(x_{0})] \leq \epsilon_{0}$ 
then $\exists \delta$, independent of $h$ such that $\forall \tau, \mbox{  } t_0 < \tau \leq t_0 + \delta R^2$, 
\[
\sup_{\stackrel{h}{\tau \leq t  \leq t_{0}+\delta \overline{R}^{2}}}  \|D^{2}u^{h}(t)\|_{L_{h}^{2}(B_{\frac{\overline{R}}{8}}(x_{0}))} \leq C
\]

\noindent \underline{Notation}: $E^{h}[u^{h}(t); B_{R}(x_{0})]$ denotes the energy of $u^{h}(t)$ indexing over $j = (j_{1}, j_{2})$ such that
$jh \in B_{R}(x_{0})$. 
\end{lemma}
\begin{proof}
The goal is to achieve uniform $H_{h}^{2}$ bounds in a reduced parabolic cylinder.  We will pass through the following intermediate stages to achieve this end:

\begin{itemize}
\item[(i)]  \emph{Show a local energy inequality}
\item[(ii)] \emph{Show $u^{h} \in L^{2}(H_{h}^{2})$ in a reduced parabolic cylinder}
\item[(iii)] \emph{Show $D^{1} u^{h} \in L^{r}(\Lh{4})$  for $r < \infty$ in a reduced parabolic cylinder}
\item[(iv)] \emph{Show $u^{h} \in L^{\infty}({H_{h}^{2}})$  in a reduced parabolic cylinder}
\end{itemize}
 Proceeding, 
\begin{itemize}
\item[(i)]
For any $\epsilon > 0$,  $E^{h}[u^{h}(t_{0}); B_{2R}(x_{0})] \leq \epsilon/2$, $\exists \delta$ such 
that $\forall t$,  $t_{0} \leq t \leq t_{0} + \delta R^{2}$, 
\begin{equation}
\label{localenergy}
\sup_{h} E^{h}[u^{h}(t); B_{R}(x_{0})] + \|\partial_t u^h\|_{L^2([t_0, t](\Lh{2}(B_{R}(x_0)))}^2 < \epsilon.
\end{equation}

\begin{subproof}
Denote by $\zeta$ the cut-off function that equals  
1 for $jh \in B_{R}(x_{0})$ and 0 for indices outside of $B_{2R}(x_{0})$  with $|D^{1} \zeta| \leq \frac{k}{R}$.  
Multiplying $({\ref{deq}})$ by $\frac{\partial u_{j}^{h}}{\partial t}\zeta^{2}$, and
summing  over $j$, 
\begin{eqnarray*}
\lefteqn{h^{2}\sum_{j}|\frac{\partial u_{j}^{h}}{\partial t}|^{2}\zeta^{2}}\\
& = & h^{2} \sum_{j} \Delta u_{j}^{h}\frac{\partial u_{j}^{h}}{\partial t}\zeta^{2} \\
& = & h^{2} \sum_{j}D_{+i}D_{-i} u_{j}^{h}\frac{\partial u_{j}^{h}}{\partial t}\zeta^{2} \\
& = & (D_{+i}D_{-i} u^{h}, \frac{\partial u_{j}^{h}}{\partial t}\zeta^{2})_{L_{h}^{2}}^{2} \\
& = &-(D_{-i} u^{h}, D_{-i}(\frac{\partial u_{j}^{h}}{\partial t}\zeta^{2}))_{L_{h}^{2}}^{2} \\
& = & h^{2} \sum_{j}\{-\zeta^{2}D_{-i} u_{j}^{h}D_{-i}\frac{\partial u_{j}^{h}}{\partial t}
    - D_{-i}\zeta^{2}D_{-i} u^{h}\frac{\partial u_{j}^{h}}{\partial t} \} \\
& \leq & h^{2} \sum_{j}-\frac{1}{2}\frac{\partial (\zeta^{2}(D_{-i}u_{j}^{h})^{2})}{\partial t} 
    +h^{2} \sum_{j} |\zeta ||D_{-i}\zeta| |D_{-i} u_{j}^{h}||\frac{\partial u_{j}^{h}}{\partial t}|\\
& \leq & h^{2} \sum_{j}-\frac{1}{2}\frac{\partial (\zeta^{2}(D_{-i}u_{j}^{h})^{2})}{\partial t}
     +\frac{1}{2}h^{2} \sum_{j} |\frac{\partial u_{j}^{h}}{\partial t}|^{2}\zeta^{2} 
     +\frac{1}{2}h^{2} \sum_{j} |D_{-i}\zeta|^{2} |D_{-i} u_{j}^{h}|^{2} \\
& \leq  & h^{2} \sum_{j} -\frac{1}{2}\frac{\partial (\zeta^{2}(D_{-i}u_{j}^{h})^{2})}{\partial t}
     +\frac{1}{2} h^{2} \sum_{j}|\frac{\partial u_{j}^{h}}{\partial t}|^{2}\zeta^{2}
     +\frac{1}{2R^{2}} h^{2} \sum_{j}  |D_{-i} u_{j}^{h}|^{2}.
\end{eqnarray*}
Since the same expression holds with $D_{+i}$ replacing $D_{-i}$, we make use of the total  energy inequality to get
\[
\frac{\partial}{\partial t}(\zeta^{2}E^{h}[u^{h}(t)])  + h^{2}\sum_{j}|\frac{\partial u_{j}^{h}}{\partial t}|^{2}\zeta^{2} \leq \frac{1}{R^{2}} E^{h}[u^{h}(t_{0})] \leq \frac{1}{R^{2}} E^{h}[f^{h}]  
\]
Integrating from $t_{0}$ to $t$,
\[
E^{h}[u^{h}(t); B_{R}(x_{0})] + \int_{t_0}^t \|\partial_t u^h \zeta\|_{\Lh{2}}^2 \leq E^{h}[u^{h}(t_{0}); B_{2R}(x_{0})] + C\frac{t-t_{0}}{R^{2}}E^{h}[f].
\] 
In particular, choosing $\delta = \frac{C\epsilon}{2E^{h}[f^{h}]}$, we get the desired result.
\end{subproof}

\item[(ii)]
\noindent \fbox{$L^{2}(H^{2})$} For $\epsilon> 0$,  $\exists \delta$ such that 
\begin{equation}
 \label{l2h2}
E^{h}[u^{h}(t); B_{2R}(x_{0})] < \epsilon \Rightarrow 
\sup_{h} \|D^{2}u^{h}\|_{L^{2}([t_{0}, t_{0} + \delta R^{2}], \Lh{2}(B_{R}(x_{0}))} \leq  
      C( \epsilon )
\end{equation}
 
 \begin{subproof}
To get bounds on terms involved in $\|D^{2}u^{h}\|_{\Lh{2}}$, it suffices to bound $\|\Delta^{h} u^{h}\|_{\Lh{2}}$ since the discrete Laplacian controls all second difference derivatives.  
For any cutoff function $\zeta$ and $\uh{j}$  solving  $({\ref{deq}})$,  $\uh{j}\zeta_{j}$ satisfies the following equation:
\begin{eqnarray}
\lefteqn{\Delta^{h}(\uh{j}\zeta_{j})}\\ & = & \{ \frac{1}{4}(\zeta_{j_{i}+1} + 2\zeta_{j_{i}} + \zeta_{j_{i}-1})\Delta^{h}\uh{j} \nonumber + 2(D_{+i}\zeta_{j} + D_{-i}\zeta_{j})(D_{+i}\uh{j} + D_{-i}\uh{j}) \nonumber\\
&  & \mbox{} +\frac{1}{4}(\uh{j_{i}+1} + 2\uh{j_{i}} + \uh{j_{i}-1})\Delta^{h}\zeta_{j} \} \\
& = & \{\frac{1}{4}(\zeta_{j_{i}+1} + 2\zeta_{j_{i}} + \zeta_{j_{i}-1})(\partial_{t}\uh{j}-\lambdah{j} \nu_{j})  + \nonumber\\ 
&  & \mbox{} 
 + 2(D_{+i}\zeta_{j} + D_{-i}\zeta_{j})(D_{+i}\uh{j} + D_{-i}\uh{j})+\frac{1}{4}(\uh{j_{i}+1} + 2\uh{j_{i}} + \uh{j_{i}-1})\Delta^{h}\zeta_{j}\} \nonumber \\
 &= &I + II + III + IV \nonumber
\end{eqnarray}
This is just the discrete analog of the expression
\[
\Delta(u\zeta) = \Delta u \zeta + \nabla u \nabla \zeta + u \Delta \zeta = (\partial_{t}u - \nabla u \cdot \nabla \nu (u))\zeta + 2\nabla u \nabla \zeta + u \Delta \zeta. 
\]
Choose $\zeta$ such that it equals 1 at indices $j$ such that $jh \in B_{R}(x_{0})$ and 0 outside of $B_{2R}(x_{0})$  with $ |D^{1} \zeta| \leq \frac{C}{R}$ and  $|\Delta^{h} \zeta| \leq \frac{C}{R^{2}}$.  Then 
\begin{itemize}
\item $\displaystyle \|I\|_{\Lh{2}(\R{2})} \leq C\|\partial_t u^h(t) \zeta\|_{\Lh{2}(\R{2})}$
\item $\displaystyle \|II\|_{\Lh{2}(\R{2})} \leq C\||D^1 u^h(t)|^2 \zeta\|_{\Lh{2}(\R{2})}$
\item $\displaystyle \|III\|_{\Lh{2}(\R{2})} \leq  C\||D^1 \zeta||D^1 u^h(t)|\|_{\Lh{2}(\R{2})} \leq \frac{C}{R} \|D^1 u^h(t)\|_{\Lh{2}(B_{2R}(x_0)} \leq \frac{C(\epsilon)}{R}$
\item $\displaystyle \|IV\|_{\Lh{2}(\R{2})} \leq  C\||\Delta^h \zeta||u^h(t)|\|_{\Lh{2}(\R{2})} \leq \frac{C}{R^2}$
\end{itemize}
where $III$ is bounded by using the local energy inequality that we attained in (i).   Now integrating on an interval of length $\delta R^2$, we can see that $III$ and $IV$ are lower order terms in the sense that 
\begin{itemize}
\item $\displaystyle \|III\|_{L^2([t_0, t_0 + \delta R^2](\Lh{2}(\R{2}))} \leq  C(\epsilon)\delta R $
\item $\displaystyle \|IV\|_{L^2([t_0, t_0 + \delta R^2](\Lh{2}(\R{2}))} \leq C\delta$
\end{itemize}
We plan to use the local space-time bound on $\partial_t u^h$ from (i) to handle $I$.  $II$, however, poses
an obstacle since there is a `missing' $\zeta$ term which prevents the step that we used in the small energy case:
\[
\||D^1 u^h|^2\|_{\Lh{2}} \leq \|D^1 u^h \|_{\Lh{4}}^2
\]
The extra work that is needed is performed in Proposition $\ref{locsob-int}$ corresponding to $s = 1, p = q = 2$.
Getting estimates on $D^{2}(u^{h}\zeta)$, and denoting the lower order terms,  $l.o.t. =   \|III\|_{\Lh{2}(\R{2})} + \|IV\|_{\Lh{2}(\R{2})}$, we have
\begin{eqnarray*}
\lefteqn{\|D^2 u^{h}(t)\zeta\|_{\Lh{2}(\R{2})}} \\
    & \leq & C\|\Delta(u^{h}(t)\zeta)\|_{\Lh{2}(\R{2})} \\
    & \leq & C\||D^{1}u^{h}(t)|^{2}\zeta\|_{\Lh{2}(\R{2})} + \|\partial_t u^h(t) \zeta\|_{\Lh{2}(\R{2})} +  l.o.t.\\
    & \leq & C\|D^{1}u^{h}(t)\|_{L_{h}^{2}(B_{2R}(x_{0}))}\|D^{1}(D^{1}u^{h}(t)\zeta)\|_{\Lh{2}(\R{2})} + \|\partial_t u^h(t)\zeta\|_{\Lh{2}(\R{2})} +  l.o.t.\\
    & \leq &C(\epsilon) \|D^{2}(u^{h}(t) \zeta)\|_{\Lh{2}(\R{2})} + \|\partial_t u^h(t) \zeta\|_{\Lh{2}(B_{2R}(x_{0}))}+  l.o.t. .
\end{eqnarray*}
For $\epsilon$ small enough, we can absorb the $\|D^{2}(u^{h}\zeta)\|_{\Lh{2}(\R{2})}$ term on the
right hand side of the equation into the left to get
\[
\|D^{2}u^{h}\|_{L_{h}^{2}(B_{R}(x_{0}))} \leq C\|\partial_t u^h(t) \zeta\|_{\Lh{2}(\R{2})}+  l.o.t.
\]
We can now take $L^{2}$ in time on the time interval $[t_{0}, t_{0}+\delta R^{2}]$ and use (i) to get
\[
\|D^{2}u^{h}\|_{L^{2}([t_{0}, t_{0}+\delta R^{2}], L_{h}^{2}(B_{R}(x_{0})))}
        \leq C(\epsilon).
\]  
\end{subproof}

\item[(iii)]
\noindent \fbox{$L^{r}(W_{h}^{1,4})$}  For $\epsilon> 0$,  $\exists \delta$ and $\tau \in (t_{0}, t_{0} + \delta R^{2}]$ such that for  $E^{h}[u^{h}(t); B_{4R}(x_{0})]  <  \epsilon$ and  $r > 1$ , 
\begin{equation}
\label{lrl4}
 \sup_{h} \|D^{1}u^{h}\|_{L^{r}[\tau, t_{0} + \delta R^{2}], \Lh{4}(B_{R}(x_{0}))}  \leq  
      C(\|u^{h}(\tau)\|_{H_{h}^{2}},  E^{h}[f^{h}]).
\end{equation}

\begin{subproof}
To get the bounds that we want for any $r > 1$, we cannot hope to start our time interval at $t_{0}$ since
all we are assuming is that $u^h(t_0) \in H_h^1$.  However, from (ii), we can find $\delta_{0}$ such that 
\[
\|D^{2}u^{h}\|_{L^{2}([t_{0}, t_{0} + \delta_{0} R^{2}], \Lh{2}(B_{2R}(x_{0}))} \leq  C.
\]
For any $t_{1} \in (t_{0}, t_{0} + \delta_{0} R^{2})$, then,  there exists $\tau \in (t_{0}, t_{1})$  such that
\[
\|D^{2}u^{h}(\tau)\|_{\Lh{2}(B_{2R}(x_{0}))}^2 \leq \frac{C}{(t_{1}-t_{0})}
\]
We can now start at $t = \tau$ with $H_{h}^{2}$ initial data.

\noindent For any choice of $\zeta$, $u_{j}^{h}$  solving  $({\ref{deq}})$,  $\uh{j}\zeta_{j}$ solves the equation
\begin{eqnarray}
\label{dcuteq}
\lefteqn{\partial_{t}(u_{j}^{h}\zeta_{j}) - \Delta^{h}(u_{j}^{h}\zeta_{j})}\\
& = & \{\frac{1}{4}(\zeta_{j_{i}+1} + 2\zeta_{j_{i}} + \zeta_{j_{i}-1})(\lambdah{j}\nu_{j})  - 2(D_{+i}\zeta_{j} + D_{-i}\zeta_{j})(D_{+i}u_{j}^{h} + D_{-i}u_{j}^{h}) \nonumber \\
&  & \mbox{} -\frac{1}{4}(u_{j_{i}+1}^{h} + 2u_{j_{i}}^{h} + u_{j_{i}-1}^{h})\Delta^{h}\zeta_{j}\nonumber \} 
\end{eqnarray}
Now we choose $\zeta$ to be the cutoff function, $1$ in a ball of radius $R$ supported in $B_{2R}(x_{0})$.  By Duhamel's,  we have
\begin{eqnarray*}
\lefteqn{\|D^{1}(u^{h}(t)\zeta)\|_{\Lh{4}(\R{2})}} \\
&\leq &  \|\Phi^h(t)(D^{1}(u^{h}(\tau)\zeta))\|_{L_{h}^{4}} + \int_{\tau}^{t} \frac{1}{(t-s)^{1/2}}\||D^{1}u^{h}(s)|^{2}\zeta\|_{\Lh{4}(\R{2})}  ds + \\
&  & \mbox{  } + \||D^1 \zeta||D^1 u^h|\|_{\Lh{4}(\R{2})} + \||u^h||\Delta^h \zeta|\|_{\Lh{4}(\R{2})} ds
\end{eqnarray*}
Using Proposition $\ref{locsob-int}$ for  $s = 2, p = 4, q = 2$, we get
\begin{eqnarray*}
\lefteqn{\|D^{1}(u^{h}(t)\zeta)\|_{L_{h}^{4}}}\\ 
& \leq &  \|\Phi^h(t-\tau)(D^{1}(u^{h}(\tau)\zeta))\|_{L_{h}^{4}} +  \int_{\tau}^{t} \frac{1}{(t-s)^{1/2}}\||D^{1}u^{h}(s)|^{2}\zeta\|_{L_{h}^{4}(\R{2})}  ds +   \\
& & \mbox{} C(R, E^h[f^h])\\
 & \leq & \|\Phi^h(t-\tau)(D^{1}(u^{h}(\tau)\zeta))\|_{L_{h}^{4}} +  C\int_{\tau}^{t} \frac{1}{(t-s)^{1/2}}
\|D^{1}(u^{h}(s)\zeta)\|_{L_{h}^{4}} \\
& & \mbox{} \cdot \|D^{2}u^{h}\|_{L_{h}^{2}(B_{2R}(x_{0}))} ds  + C(R, E^h[f^h])
\end{eqnarray*}  
Now taking $L^{r}$ in time of both sides, we know by Young's Inequality that we can take $\|D^{2}u^{h}\|_{L_{h}^{2}(B_{2R}(x_{0}))}$ in $L^{2}$ and $\|D^{1}(u^{h}(s)\zeta)\|_{L_{h}^{4}}$ in $L^{r}$;  the first term is  bounded by $C(\epsilon)$ by (ii) and the second is exactly the quantity that we are trying to bound.  This entire product can be absorbed into the left hand side.  The linear term  $\|\Phi^h(t)(D^{1}(u^{h}(\tau)\zeta))\|_{L_{h}^{4}}$ is in $L^{r}$ since $u^{h}$ is $H_{h}^{2}$ at $t = \tau$.  Therefore, 
\[
\|D^{1}(u^{h}\zeta)\|_{L^{r}([\tau, t_{0}+\delta R^{2}], L_{h}^{4}(\R{2})} \leq  C(\|u^{h}(\tau)\|_{H_{h}^{2}}, E^h[f^h]).
\]

\end{subproof}

\item[(iv)]
\noindent \fbox{$L^{\infty}(H_{h}^{2})$} 
For $\epsilon> 0$,  $\exists \delta$ and $\tau \in (t_{0}, t_{0}+\delta R^{2})$ such that for $s < \infty$, 
\begin{equation}
\label{linftyh2}
E^{h}[u^{h}(t); B_{8R}(x_{0})] < \epsilon \Rightarrow 
\sup_{\stackrel{h}{\tau \leq t \leq t_{0} + \delta R^{2}}} \|D^{2}u^{h}\|_{\Lh{2}(B_{R}(x_{0}))} \leq  
      C(\|u^{h}(\tau)\|_{H_{h}^{2}},  E^{h}[f^{h}])
\end{equation}

\begin{subproof}
As in (iii), we can use the $L^{2}(H_{h}^{2})$ bounds to find $\tau$ such that $u^{h}(\tau) \in H_{h}^{2}$.   For second derivative estimates on $u^{h}\zeta$, we need to differentiate $({\ref{dcuteq}})$ and use the
estimates $({\ref{de2}})$ by taking $p = 4, q = p' = 4/3$.    Using Observation 1 and 2 on $\lambdah{j}$ and the smoothness of $\nu$, 
\[
D^{1}(\lambdah{j} \nu_{j})
 \letilde   |D^{1}\lambdah{j}| + |\lambdah{j}||D^{1}\nu_{j}|  
 \letilde  |D^{1}u_{l}|^{3} +  |D^{1}u_{l}||D^{2}u_{l}|
\]
To avoid cluttering the basic idea with overwhelming technical
detail, we highlight the crucial term:
\begin{eqnarray*}
\lefteqn{\|D^{2}(u^{h}(t)\zeta)\|_{\Lh{2}(\R{2})}} \\
  & \leq & \|\Phi^h(t-\tau)D^{2}(u^{h}(\tau)\zeta)\|_{\Lh{2}(\R{2})}+ \\
 & & \mbox{  } + \int_{\tau}^{t} \frac{1}{(t-s)^{3/4}}\|\zeta D^{1}(\lambda^{h}(s)\nu(s))|_{\Lh{4/3}(\R{2})}ds + C(R, E^h[f^h]) \\
  & \leq &  C\|D^{2}(u^{h}(\tau)\zeta)\|_{\Lh{2}(\R{2})}+ C\int_{\tau}^{t} \frac{1}{(t-s)^{3/4}}
 \{ \||D^{1}u^{h}(s)|^{3}\|_{\Lh{4/3}(\R{2})} +\\
 & & \mbox{}   ||D^{2}u^{h}(s)||D^{1}u^{h}\zeta)\|_{\Lh{4/3}(\R{2})} ds\} + C(R, E^h[f^h]) \\
  & \leq &  C\|D^{2}(u^{h}(\tau)\zeta)\|_{\Lh{2}(\R{2})}+C \int_{\tau}^{t} \frac{1}{(t-s)^{3/4}}
  \{\|D^{1}u^{h}(s)\|_{\Lh{2}}\||D^{1}u^{h}|^{2}\zeta\|_{\Lh{4}} +  \\ & \mbox{  } &
  \|D^{1}u^{h}(s)\|_{L_{h}^{4}(B_{2R}(x_{0}))}\|D^{2}(u^{h}(s)\zeta)\|_{\Lh{2}(\R{2})}\}  ds + C(R, E^h[f^h])\\
  & \leq &  C\|D^{2}(u^{h}(\tau)\zeta)\|_{\Lh{2}(\R{2})}+ C\int_{\tau}^{t} \frac{1}{(t-s)^{3/4}}
  \|D^{1}u^{h}(s)\|_{L_{h}^{4}(B_{2R}(x_{0}))} \\
  & & \mbox{ } \cdot \|D^{2}(u^{h}(s)\zeta)\|_{\Lh{2}(\R{2})}  ds+C(R, E^h[f^h])
\end{eqnarray*}
Now we can take $L^{\infty}$ in time of this expression by taking $\|D^{2}(u^{h}(s)\zeta)\|_{\Lh{2}(\R{2})}$ in $L^{\infty}$ and $\|D^{1}u^{h}(s)\|_{L_{h}^{4}(B_{2R}(x_{0}))}$ in $L^{r}$ for $ r > 4$.  From this we can conclude $L^{\infty}(H_{h}^{2})$ bounds in the parabolic cylinder $B_{R}(x_{0}) \times [\tau, t_{0} + \delta R^{2}]$.  
\end{subproof}
\end{itemize}
\end{proof} 

\begin{proof} (\emph{Theorem \ref{fehhpthm1}})
The construction of a solution to $(\ref{ceq})$ begins, as before, at the level of the discrete approximations.  Existence for any fixed $h$ is not an issue.  Moreover, since our initial data is of finite energy, we know from the energy inequality that any solution remains of finite energy.  However, this energy may concentrate which in turn invalidates an immediate argument for the convergence of the interpolants.  As seen in the case of small energy initial data, the convergence argument depended on 
having control on the second derivatives, or uniform control on the gradient.  What Lemma $\ref{fehhpthm2}$ provides is local control of these quantities.  The key extra step that is needed
is an estimate of the size of the energy concentration set.    From this we can conclude that the limit of the $H^{1}$ interpolants can be extended across the concentration sets in such
a way that it is a partially regular solution of the harmonic map heat flow problem. 

\subsubsection*{Step 1: Attaining Bounds} This is the content of Lemma $\ref{fehhpthm2}$.

\subsubsection*{Step 2: Estimating the Size of the Concentration Set}
By invariance of the equation $\ref{deq}$ by time translation, it suffices to estimate the concentration 
set on the interval $[0,1]$.  Fix $\epsilon_{0}$ and $\delta$ of Theorem $\ref{fehhpthm1}$.  Since $\delta$ is independent of $h$, it makes sense, for fixed $R_0$, to
partition $[0,1]$ into intervals $I_0, I_1, \ldots, I_{2/\delta R_0^2}$ of length $\delta R_0^2/2$.
Denote  $T_{j}$ to be the base time of $I_{j}$.
Define $\beta: [0,1] \rightarrow \{0, \ldots,  2/\delta R_0^2\}$, assigning to every $t$ in $[0,1]$ the index $j$ such that 
$t \in I_j$.      For arbitrary $R$ and  any point
$(x_{0}, t_{0}) \in \R{2}\times (0, 1]$ denote by $P_{R}(x_{0}, t_{0})$ the parabolic cylinder 
\[
P_{R}(x_{0}, t_{0}) = \{(x, t) | x \in B_{R}(x_{0}), t \in [t_0-\frac{\delta R^2}{2}, t_0+\frac{\delta R^2}{2}]\}
\]
\noindent Define the singular set
\[
\Sigma = \bigcap_{R> 0} \{(x, t) \in \R{2}\times \Rplus : \liminf_{h\rightarrow 0} E^{h}[u^{h}, P_{R}(x, t)]  >  \epsilon_{0}\}
\]
Then 
$(x_k, t_k) \in \Sigma$ implies, by Lemma  \ref{fehhpthm2}, 
that for every $R \leq R_0$ and every $h$ sufficiently small,
\[
E^{h}[u^{h}(t_k - \frac{\delta R^2}{2}), B_{2R}(x_{k})] > \epsilon_{0}/2.
\]  
For any $(x_k, t_k)$, associate $R_k \leq R_0$  satisfying
\[
t_k - \frac{\delta R_k^2}{2} = T_{\beta(t_k)} 
\]
Then 
\[
E^h[u^h(T_{\beta(t_k)}), B_{2 R_k(x_k)}] >  \epsilon_0/2.
\]
Now these concentration balls  have been `snapped to grid' associated to $R_0$.  
At any time slice,  a direct consequence of the energy bound is that there is a finite number of balls of 
radius $2R$ in which the energy exceeds $\epsilon_{0}/2$, independent of $R$.   Specifically, there
can be no more than $E^{h}[f^{h}]/(2\epsilon_{0})$ such balls at any time slice for any $R$. 
The balls $B_{2R_k}$ at time $T_j$ cover the bases of the concentration cylinders in $I_j$ in the sense of Vitali, 
so by Vitali's covering lemma,  we can take these balls to be disjoint.   Now adding up the contribution of all the concentration cylinders over the $2/(\delta R_0^{2})$ intervals, we have 
\[
\Sigma_{k} Vol(C_{k}) \leq \frac{1}{2}\delta R_0^2 \Sigma_{k} R_{k}^{2} \leq  \frac{1}{2}\delta R_0^{2} \frac{E^{h}[f^{h}]}{2\epsilon_{0}}\delta R_0^{2} \frac{2}{\delta R_0^2}
 \]
Thus
\[
\Sigma_{k} R_{k}^{2} \letilde  \frac{E^{h}[f^{h}]}{2\delta \epsilon_{0}}.
\]
We can conclude that the concentration set $\Sigma$ has locally finite $2$-dimensional Hausdorff measure with respect to the parabolic metric.

\subsubsection*{Step 3: Convergence of the Interpolants}
For $(x_{0}, t_{0}) \in \Sigma^{c}$, there exists $R > 0$ such that for some unbounded sequence $h$, 
\[
E^{h}[u^{h}, P_{R}(x_{0}, t_{0})]  \leq  \epsilon_{0}.
\]
Equivalently, 
\[
E^{h}[u^{h}(T_{t}), B_{2R}(x_{0})] \leq \epsilon_{0}/2,
\]  
so by Lemma $\ref{fehhpthm1}$, $\|D^{2}u^{h}\|_{\Lh{2}} \leq C$ in a uniform neighborhood Q of $(x_0, t_0)$. 
Associate with each sequence $\uh{j}$ the $H^{1}$ interpolant $p^{h}$.  From Lemma $\ref{phlemma}$, 
$p^{h} \rightarrow p $ in $H^{1}(Q)$.   By the argument given in the
global small energy case we have, for $\zeta \in C_{0}^{\infty}(Q)$, that
\[
\int_{Q} \partial_{t}p\cdot \zeta + \nabla p \cdot \nabla \zeta - \nabla p \cdot \nabla\nu(p))\nu(p) \cdot \zeta
 = 0. 
 \] 
Off $\Sigma$, then, $p$ solves ($\ref{ceq}$).  

To show that $p$ solves the equation on $\Rplus\times \R{2} $, fix a compact set $Q \subset \subset \Rplus \times \R{2}$.  As in the covering argument above, we can, for fixed $R$,  partition up the time interval occupied by $Q$ into intervals of length $\delta R^{2}$.  Without loss of generality, we take the total time interval to be of unit length.   Given $R > 0$, let $\{Q_{i,l} = P_{R}(x_i,l, t_l)\}$ be a covering of 
$\Sigma \cap Q$ by parabolic cylinders.  Denote by $N_l$ the number of parabolic cylinders in each time interval $I_l$.  By the energy inequality, $N_l$ is bounded by some finite number $N_{0}$ independent of $R$, so the number of such cylinders is less than $ N_{0}/(\delta R^{2})$.   

Consider  smooth cut-off functions $0 \leq \eta_{1}, \eta_{2} \leq 1$ identically 0 in $B_{1}(0)$ and 
$1$ outside of $B_{2}(0)$.  Define the cutoff function
\[
\eta_R =  \sum_{l = 1}^{1/\delta R^{2}}\sum_{i =1}^{N_{l}} \eta_{1}(\frac{x-x_{i,l}}{R})\eta_{2}(\frac{t-t_{l}}{\delta R^{2}}).
\]
Let $\phi \in C_{0}^{\infty}(Q)$. Then setting $\zeta = \phi \eta_{R}$, 
\[
0 = \int_{Q} \{\partial_{t}p\cdot \phi + \nabla p \cdot \nabla \phi - \nabla p \cdot \nabla\nu(p))\nu(p) \cdot \phi\}\eta_{R}+ L
\]
where $L$  can be bounded by 
\begin{eqnarray*}
|L| &\leq & \int_{Q} |\nabla p ||\phi||\nabla \eta_{R}|\\
    &\leq & \int_{Q} |\nabla p ||\nabla_{x}(\displaystyle \sum_{l = 1}^{1/\delta R^2}
   \sum_{i = 1}^{N_{l}} \eta_{1}(\frac{x-x_{i,l}}{R})\eta_{2}(\frac{t-t_{l}}{\delta R^{2}}))| dx dt \\
    & \leq & C \sum_{l = 1}^{1/\delta R^2}
   \sum_{i = 1}^{N_{l}} \frac{1}{R} \int_{Q_{i,l}} |\nabla p|  dx dt\\
   & \leq &  C\sum_{l = 1}^{1/\delta R^2}
   \sum_{i = 1}^{N_{l}} \frac{1}{R} (\int_{Q_{i,l}} |\nabla p|^{2})^{1/2} (\int_{Q_{i,l}} dx dt)^{1/2}  
\end{eqnarray*}
\begin{eqnarray*}
   & \leq &   C\sum_{l = 1}^{1/\delta R^2}\sum_{i = 1}^{N_l} R (\int_{Q_{i,l}} |\nabla p|^{2})^{1/2}  \\
   & \leq &  C(\sum_{i = 1}^{N_l} \sum_{l = 1}^{1/\delta R^2} R^2)^{1/2}(\sum_{i = 1}^{N_0} \sum_{l = 1}^{1/\delta R^2} \int_{Q_{i,l}} |\nabla p|^{2})^{1/2} \\
   & \leq & CN_0^{1/2} (\sum_{i = 1}^{N_0} \sum_{l = 1}^{1/\delta R^2} \int_{Q_{i,l}} |\nabla p|^{2})^{1/2}
\end{eqnarray*}
Taking $R \rightarrow 0$, $|L| \rightarrow 0$  by absolute continuity so $p$ weakly solves $({\ref{ceq}})$ in $\Rplus \times \R{2}$. 
\end{proof}

\section{LLG}
\label{llg}
We proceed to prove our main result,  Theorem \ref{fellgthm1}.   The construction of a solution to this problem will follow the same path as that for the harmonic heat flow.  The semi-discrete approximation for LLG is given by  (\ref{dllg}).  
The methods of the previous section heavily depended on estimates for the linear equation.   One additional ingredient  is  needed for the LLG:   a reformulation of   (\ref{dllg}) that yields a system that is linear in its highest order term.   In the continuous setting, this reformulation is achieved by  showing that the coordinates of derivatives of $u$ solving ($\ref{llgequation}$) for $\alpha = 0$ onto a moving frame satisfy a system whose highest order term is a linear Schr\"odinger term.   The concept of frames is reviewed in Appendix B.  For our (discrete) setting, we achieve such a formulation by  replacing derivatives of $u$ with projections onto the tangent space of difference derivatives.  
\subsection{The Linearized Discrete System}
\label{linearizedsystem}
\subsubsection*{Notation}
Difference derivative lose their tangency, so we separate out the tangent and the orthogonal portions.  Denote
\[
\begin{array}{ccc}
D_{+k}\uh{j} = \vplusk{j} + \wplusk{j},  &  & D_{-k}\uh{j} = \vmink{j} + \wmink{j}
\end{array}
\]
where $\vplusk{j}, \vmink{j} \in  \TN$ and $\wplusk{j}, \wmink{j} \perp \TN$ .  
The tangent portion can be further decomposed as follows. We can assume that
$\mathcal{N}$ admits a global frame $\{e^1, e^2 = \nu \wedge e^1\}$.  Justification for this is provided In Appendix B.  Then $e_j^i = e^i (\uh{j})$.  
Introduce the complex notation 
\[
\begin{array}{ccccc}
 i: T\mathcal{N} \rightarrow T\mathcal{N},  & ie^1 = e^2, &  &   ie^2 = -ie^1. 
\end{array}
\] 
Letting $e = e^1$ and $\qplusk{j} = \qpluski{j}{1} + i \qpluski{j}{2}$, $\qmink{j} = \qminki{j}{1} + i\qminki{j}{2}$,
\[
 \vplusk{j} =  \PmN(D_{+k}u_{j})=  \qplusk{j}e_{j}, \mbox{     }  \vmink{j} = \PmN(D_{-k}u_{j}) =  \qmink{j}e_{j}
\]
where \PmN denotes the orthogonal projection onto $T_{u_j}\mathcal{N}$.  
The orthogonal portion can be written as
\[
 \wplusk{j} =    \aplusk{j}\nu_{j}, \mbox{     }  \vmink{j}  =  \amink{j}\nu_{j}
\]
$\mN$  is a smooth hypersurface so the unit normal $\nu$ is a smooth function.  Denote 
by 
\[
c_{\nu} = \|\nabla \nu\|_{L^{\infty}}.
\]  
 Since $\mN$ is smooth, $e$ is smooth so we denote 
$c_{e} = \|\nabla e\|_{L^{\infty}}.$
\subsubsection*{Observations}
\begin{itemize}
\item[(1)] $|\vplusk{j}| = |\qplusk{j}|$ and $|\vmink{j}| = |\qmink{j}|$
\item[(2)]  $\aplusk{j}= D_{+k}\uh{j}\cdot \nu_{j} = O(h|D_{+k}\uh{j}|^{2}).$  Similarly,  $\amink{j}=  O(h|D_{-k}\uh{j}|^{2}).$ 
\item[(3)] $\lambdah{j} = O(|D^{1}\uh{j}|^{2}).$ 
\item[(4)] For suitable $\epsilon$, and $|\uh{j_{k}+1}-\uh{j}|^2 < \epsilon$, 
\[
\begin{array}{cc}
|D_{+k}\uh{j}|^2 = O(|\vplusk{j}|^2),   & D_{+k}\uh{j}|^2 = O(|\vplusk{j}|^2).
\end{array}
\]
This follows from the following manipulation:
\begin{eqnarray*}
|D_{+k}\uh{j}|^{2} \leq 2(|\vplus{j}|^{2} + |\wplus{j}|^{2}) & = &(|\vplus{j}|^{2} +  |\aplus{j}u_{j}|^{2}) \\
&\leq & 2(|\vplus{j}|^{2}  + C||hD_{+k}\uh{j}|^{2}) \\
&\leq & 2(|\vplus{j}|^{2}  + C||\uh{j_{k}+1}-\uh{j}|^2 |D_{+k}\uh{j}|^{2}),
\end{eqnarray*}
where we have used (1) in the last line.
Using a similar argument, for suitable $\epsilon$,  

\noindent $|\uh{j_{k}+1}-\uh{j}|^2 < \epsilon \Rightarrow$
\[
\begin{array}{cc}
|D_{+k}\uh{j}|^{2} = O(|\qplusk{j}|^{2}), &  |D_{-k}\uh{j}|^{2} = O(|\qmink{j}|^{2}) \\
\aplusk{j}= O(|\qplusk{j}|), &   \amink{j}= O(|\qmink{j}|)
\end{array}
\]
\item[(5)] Since $\nu_{j} =  \nu(\uh{j})$
 \[
 \begin{array}{cc}
\partial_{t}\nu_{j} = \nabla_{u}\nu \partial_{t}\uh{j} = O(|\partial_{t}\uh{j}|), &  D_{k}\nu_{j} = O(|D_{k}\uh{j}|)
\end{array}
\]
\item[(6)]  $e$ is smooth and $e_j = e(u_j)$ so 
\[
\begin{array}{cc}
\partial_{t}e_{j} = \nabla_{u}e \partial_{t}\uh{j} = O(|\partial_{t}\uh{j}|), & D_{k}e_{j}  = O(|D_{k}\uh{j}|)
\end{array}
\]
\end{itemize}
\begin{lemma}
\label{dhasimoto}
For $\epsilon > 0$, choose $h$ such that 
\[
\begin{array}{cc}
|\uh{j_{i}+1} - \uh{j}|^{2} < \epsilon,   &  | \uh{j} - \uh{j_{i}-1}|^{2} < \epsilon.
\end{array} 
\]
Then $q_{j} = \qplusk{j}$ or $\qmink{j}$  satisfies
\[
\partial_{t}q_{j}  = i\Delta^{h} q_{j} + \Delta^{h} q_{j} + F^{h}(q_{l}, D^{1}q_{l}),
\]
where $F^{h}(q, D^{1}q) = O(|q|^{3}, |q||D^{1} q|)$ and $l$ indexes over nearest neighbors of $j$. 
\end{lemma}

\begin{proof}

\noindent Let $u_j = \uh{j}$, $q_{j} = \qplusk{j}$ and $v_{j} = \vplusk{j}$. To get an equation in terms of $q_{j}$,  take $D_{+k}$ of (\ref{dllg}) and then take the inner product of the resulting expression  with $e_{j}$.   Let 
\[
I =  D_{+k} (\partial_{t}u_{j})\cdot e_{j} = \partial_t (D_{+k}u_j)\cdot e_j = II + III + IV . 
\]
where
\begin{itemize}
\item[$II$] $\displaystyle  = (\nu_{j} \wedge D_{+k} \Laph \uh{j})\cdot e_{j}$
\item[$III$] $\displaystyle =  (D_{+k}\nu_{j} \wedge D_{+k}^{2}\uh{j})\cdot e_{j}$
\item[$IV$] $\displaystyle = D_{+k} ( \Laph \uh{j} + \lambdah{j}\nu_{j})\cdot e_{j}$
\end{itemize} 

\begin{itemize}
\item[ ]
\begin{eqnarray*} 
I & = & \partial_{t}(D_{+k}u_{j})\cdot e_{j} \\
& = &  \partial_{t}(v_{j} + a_{j}\nu_{j})\cdot e_{j}  \\
& = &  ( \partial_{t}(q_{j}e_{j}) + a_{j}\partial_{t}\nu_{j})\cdot e_{j}  \\
& = & \partial_{t}q_{j} + q_{j}\partial_{t}e_{j}\cdot e_{j} +   a_{j}c_{\nu} \partial_{t}u_{j}\cdot e_{j}\\
& = &  \partial_{t}q_{j} +(c_e q_{j} + c_{\nu}a_{j}) D_{+k}u_{j}\cdot \nu_{j})\partial_{t}u_{j}\cdot  e_{j} \\
& = & \partial_{t}q_{j} + (c_e q_{j} + c_{\nu} a_{j}) D_{+k}u_{j}\cdot \nu_{j})(\nu_{j}\wedge  \Laph u_{j}\cdot e_{j}+ \Laph u_{j}\cdot e_{j})
\end{eqnarray*}
Each term except $\partial_{t}q_{j}$ belongs in the nonlinear term $F(q_{j}, D^{1}q_{j})$.  We need to show that the 
magnitude of each term is of the order claimed.
\begin{itemize}
\item[$\bullet$] $(\nu_{j} \wedge \Laph u_{j})\cdot e_{j} =   ( \nu_{j} \wedge D_{-i}v_{j})\cdot e_{j}  + ( \nu_{j} \wedge D_{-i}w_{j})\cdot e_{j} $
\begin{eqnarray*}
| \nu_{j} \wedge D_{-i}v_{j})\cdot e_{j} | 
& = &  |\nu_{j} \wedge (D_{-i}q_{j} + q_{j_{i}-1}D_{-i}e_{j})\cdot e_{j}| \\
& = &| D_{-i}q_{j}(\nu_{j}\wedge e_{j})\cdot e_{j} + c_e q_{j_{i}-1}(\nu_{j}\wedge \vmina{j})\cdot e_{j}| \\
& = & |i D_{-i}q_j + c_e q_{j_{i}-1}(\nu_{j}\wedge \vmina{j})\cdot e_{j}| \\
& = & O(|D^1 q_l|, |q_{l}|^{2})
\end{eqnarray*}
\item[ ] 
\begin{eqnarray*}
|( \nu_{j} \wedge D_{-i}w_{j})\cdot e_{j} |
& = &| \nu_{j}\wedge D_{-i}(a_{j}\nu_{j})\cdot e_{j}| \\
& = &| (\nu_{j}\wedge a_{j_{i}-1}D_{-i}\nu_{j})\cdot e_{j}|\\
& = & |c_{\nu} a_{j_{i}-1}(\nu_{j}\wedge D_{-i}u_{j})\cdot e_{j} |\\
& = & O(h|D_{-i}u_{j}|^{2}|D_{-i}u_{j}| \\
& = & O((u_{j_{i}}-u_{j_{i}-1})|D_{-k}u_{j}|^{2}|\\
& = & O(|q_{l}|^{2})
\end{eqnarray*}
\item[$\bullet$]  Similarly, \begin{eqnarray*}
|\Laph u_{j}\cdot e_{j} |&=& |D_{-i}(v_{j}+w_{j})\cdot e_{j}| \\
&=& |D_{-i}q_{j} + (c_e q_{j_{i}-1} + c_{\nu} a_{j_{i}-1})D_{-i}u_{j}\cdot e_{j} | \\
&=& O(|D^{1}q_{l}|, |q_{l}|^{2}).
\end{eqnarray*}
\end{itemize}
Multiplying both terms by $|(c_e q_{j} + c_{\nu} a_{j})| = O(|q_{j}|)$, the contribution that $I$ has to $F$ is $O(|q_{l}|^{3} , |q_{l}||D^{1}q_{l}|)$. 
\item[ ] \begin{eqnarray*}
II & = & (\nu_{j} \wedge D_{+k} \Laph u_{j})\cdot e_{j} \\
& = & (\nu_{j} \wedge \Laph(v_{j} + w_{j}))\cdot e_{j} \\
& = & (\nu_{j} \wedge \Laph v_{j})\cdot e_{j} + (\nu_{j}\wedge \Laph w_{j})\cdot e_{j}
\end{eqnarray*}
\begin{itemize} 
\item[$\bullet$] The last term belongs to the nonlinear term $F(q_{j}, D^{1}q_{j})$.  
\begin{eqnarray*}
| (\nu_{j}\wedge \Laph w_{j})\cdot e_{j}| 
& = & \frac{1}{h}|(\nu_{j} \wedge (D_{+i}-D_{-i})w_{j})\cdot e_{j}|\\
& = & \frac{1}{h}|(\nu_{j}\wedge a_{j_{i}+1}D_{+i}u_{j}-\nu_{j}\wedge a_{j_{i}-1}D_{-i}u_{j})\cdot e_{j}|\\
& = & \frac{1}{h}|(\nu_{j}\wedge a_{j_{i}+1}\vplusa{j}-\nu_{j}\wedge a_{j_{i}-1}\vmina{j})\cdot e_{j}|\\
& = & O(\frac{1}{h}h\{|D^{1}u_{j}|^{2}|\vplusa{j}, |D^{1}u_{j}|^{2}|\vmina{j}\}) \\
& = & O(|q_{l}|^{3})
\end{eqnarray*}
\item[$\bullet$] To handle the term $(\nu_{j} \wedge \Laph v_{j})\cdot e_{j}$, consider first $\Laph v_{j}$. 
\noindent \begin{eqnarray*}
 \Laph v_{j}
 & = & \Laph (q_{j}e_{j}) \\
 & = & \frac{1}{h}(D_{+i}(q_{j}e_{j}) - D_{-i}(q_{j}e_{j})) \\
 & = & \Laph q_{j}e_{j}+ \frac{1}{h}(q_{j_{i}+1}D_{+i}e_{j} - q_{j_{i}-1}D_{-i}e_{j})
\end{eqnarray*}
\noindent The first term yields the linear Schr\"odinger term
\[
(\nu_{j} \wedge \Laph q_{j} e_{j})\cdot e_{j} = \Laph q_{j}(\nu_{j} \wedge e_{j})\cdot e_{j} =  i\Laph q_{j}.
\]
\noindent The magnitude of the second term can be estimated as follows:
\begin{eqnarray*}
\lefteqn{|(\nu_{j}\wedge \frac{1}{h}(q_{j_{i}+1}D_{+i}e_{j} - q_{j_{i}-1}D_{-i}e_{j}))\cdot e_{j}| } \\
& = & |\frac{1}{h}c_e\{q_{j_{i}+1}(\nu_{j} \wedge D_{+i}u_{j}\cdot e_{j} - q_{j_{i}-1}(\nu_{j} \wedge D_{-i}u_{j})\cdot e_{j}\}| \\
& = &  |\frac{1}{h}c_e\{q_{j_{i}+1}(\nu_{j} \wedge \vplusa{j})\cdot e_{j} - q_{j_{i}-1}(\nu_{j} \wedge \vmina{j})\cdot  e_{j}\}| \\
& = &  |\frac{1}{h}c_e(q_{j_{i}+1}i\qplusa{j} - q_{j_{i}-1}i\qmina{j})| \\
& = & |c_e(q_{j_{i}+1}D_{0\i}q_{j} + \qmina{j}D_{0\i}q_{j})|\\
& = & O(|q_{l}||D^{1}q_{l}|).
\end{eqnarray*}
\end{itemize}
\item[ ]\begin{eqnarray*}
III & = & D_{+k}\nu_{j} \wedge D_{+k}^{2}u_{j}\\
& = & c_{\nu} (D_{+k}u_{j}\wedge \frac{1}{h} D_{+k}u_{j+1})\\
& = & c_{\nu}\frac{1}{h}(D_{-k}u_{j+1}\wedge D_{+k}u_{j+1})\\
& = & c_{\nu}\frac{1}{h}((\vmin{j+1} +  \wmin{j+1})\wedge  (\vplus{j+1}+ \wplus{j+1}))\\
& = & c_{\nu}\frac{1}{h}(\vmin{j+1}\wedge \wplus{j+1} + \wmin{j+1}\wedge \vplus{j+1})\\
& = & c_{\nu}\frac{1}{h}(\vmin{j+1}\wedge \wplus{j+1} + \wmin{j+1}\wedge \vplus{j+1})\\
& = & c_{\nu}\frac{1}{h}((\qmin{j+1}e_{j})\wedge (\aplus{j+1}\nu_{j+1}) + (\amin{j+1}\nu_{j+1})\wedge (\qplus{j+1}e_{j}))\\
& = & c_{\nu}\frac{1}{h}(\qmin{j+1}\aplus{j+1}(-ie_{j+1}) + \amin{j+1}\qplus{j+1}(ie_{j+1}))
\end{eqnarray*}
This term, which contributes to $F$, has the correct order since $e_{j+1} = e_{j} + O(h)$, and 
\[
|\frac{1}{h}(\qmin{j+1}\aplus{j+1}(-ie_{j+1}) + \amin{j+1}\qplus{j+1}(ie_{j+1}))\cdot e_{j+1}|  = O(|q_{l}^{3}|).
\]
\item[ ]
\[
IV = D_{+k} ( \Laph u_{j} + \lambdah{j}\nu_{j}) \cdot e_{j} = 
 \Laph v_{j}\cdot  e_{j} + \Laph w_{j}\cdot e_{j} + \lambdah{j+1}D_{+k}\nu_{j}\cdot e_{j}
\]
\begin{itemize}
\item[$\bullet$]  $\displaystyle \Laph v_{j}\cdot e_{j}  =  \Laph q_{j} + \frac{1}{h}q_{j_{i}+1}D_{+i}e_{j} - q_{j_{i}-1}D_{-i}e_{j}\cdot e_{j} $
\begin{eqnarray*}
 |\frac{1}{h}(q_{j_{i}+1}D_{+i}e_{j} - q_{j_{i}-1}D_{-i}e_{j})\cdot e_{j}|
 & = & |c_e \frac{1}{h}(q_{j_{i}+1}\vplus{j} - q_{j_{i}-1}\vmin{j})\cdot e_{j}|\\
 & = & O(|q_{l}||D^{1}q_{l}|).
\end{eqnarray*}
\item[$\bullet$] $\displaystyle \Laph w_{j}\cdot e_{j} = \frac{1}{h}(a_{j_{i}+1}D_{+i}\nu_{j}-a_{j_{i}-1}D_{-i}\nu_{j})\cdot e_{j}$
\begin{eqnarray*}
\frac{1}{h}(a_{j_{i}+1}D_{+i}\nu_{j}-a_{j_{i}-1}D_{-i}\nu_{j})\cdot e_{j}
& = & C\frac{1}{h}(a_{j_{i}+1}D_{+i}u_{j}-a_{j_{i}-1}D_{-i}u_{j})\cdot e_{j} \\
& = & C\frac{1}{h}(a_{j_{i}+1}\vplus{j}-a_{j_{i}-1}\vmin{j})\cdot e_{j} \\
& = & O(|q_{l}||D^{1}q_{l}|).
\end{eqnarray*}
\item[$\bullet$] $\displaystyle \lambdah{j+1}D_{+k}\nu_{j}\cdot e_{j}= O(|q_{l}|^{3}) $
\end{itemize}
\end{itemize}
\end{proof}
 
\subsection{Proof of Main Result}
Equipped with Lemma \ref{dhasimoto}, we can follow the steps illustrated for the harmonic map heat flow to prove our main result.   
\begin{lemma}
\label{fellgthm2}
There is an $\epsilon_{0} > 0$ such that if $E^{h}[u^{h}(t_{0}); B_{\overline{R}}(x_{0})] \leq \epsilon_{0}$ 
then $\exists \delta$, independent of $h$ such that
\[
\sup_{\stackrel{h}{t_0 < t  \leq t_{0}+\delta \overline{R}^{2}}}  \|D^{2}u^{h}(t)\|_{L_{h}^{2}(B_{\frac{\overline{R}}{8}}(x_{0}))} \leq C
\]
\end{lemma}
\begin{proof}
From  ($\ref{dgeneral2}$), the uniform bounds on energy and  $\|\partial_{t}u^{h}\|_{\Lh{2}}$ are immediate.  Multiplying $(\ref{dgeneral2})$ by  $\displaystyle \dt{\uh{j}}$, summing over $j$ and integrating in $t$, we have
\[
E^{h}[u^{h}(t)]  + \|\dt{u^{h}}\|_{L^{2}(\Lh{2})} ^{2} \leq E^{h}[f^{h}].
\]
We will achieve $ L^{\infty}(H_{h}^{2})$ bounds in the incremental way that was illustrated for the harmonic map heat flow:
\begin{itemize}
\item[(i)] A local energy inequality follows exactly as in the proof of $\ref{localenergy}$ using  ($\ref{dgeneral2}$)
\item[(ii)] \fbox{$L^{2}(H_{h}^{2})$}  Also using  ($\ref{dgeneral2}$), nothing changes in the proof of
this estimate except that $\partial_t u_j$ has to be replaced by
$ \partial_{t}\uh{j} - \nu_{j}\wedge \partial_{t}\uh{j}.$
\item[(iii)] \fbox{$L^{r}(W_{h}^{1,4})$} \hspace{0.1in}   The statement that we will show is:

\noindent For $\epsilon> 0$,  $\exists \delta$ and $\tau \in (t_{0}, t_{0} + \delta R^{2}]$ such that for $r > 1$ and $u^{h}$ solving ($\ref{dllg}$)
\[
E^{h}[u^{h}(t); B_{4R}(x_{0})] < \epsilon \Rightarrow 
\sup_{h} \|D^{1}u^{h}\|_{L^{r}[\tau, t_{0} + \delta R^{2}], \Lh{4}(B_{R}(x_{0}))} \leq  
      c(\|u^{h}(\tau)\|_{H_{h}^{2}},  E^{h}[f^{h}])
\]
\begin{subproof}
We can use the $L^{2}(H_{h}^{2})$ bounds to show that for any  $t_{1}\in (t_{0}, t_{0} + \delta R^{2})$,  there exists $\tau \in (t_{0}, t_{1})$  such that
\[
\|D^{2}u^{h}(\tau)\|_{\Lh{2}(B_{2R}(x_{0}))} \leq \frac{C}{(t_{1}-t_{0})}
\]
We can now start at $t = \tau$ with $H_{h}^{2}$ initial data.
Since the intention is to use linear estimates on the fundamental solution,  we will resort to the linearized
discrete system given by Lemma $\ref{dhasimoto}$.   Since $E^{h}[u^{h}(t); B_{4R}(x_{0})] < \epsilon$,
we have in particular that 
\[
\begin{array}{ccc}
|\uh{j_{i}+1} - \uh{j}|^{2} < \epsilon, &  &| \uh{j} - \uh{j_{i}-1}|^{2} < \epsilon. 
\end{array}
\]
so the condition of the lemma is satisfied.
\noindent Choose $\zeta$ to be the spatial cutoff function, $1$ in a ball of radius $R$ supported in $B_{2R}(x_{0})$. Then  for $q_{j} = \qplusk{j}$ or $\qmink{j}$, $q_{j}\zeta_{j}$ solves the equation
\begin{eqnarray*}
\lefteqn{\nonumber \partial_{t}(q_{j}\zeta_{j}) - (1+i)\Delta^{h}(q_{j}\zeta_{j})}\\
& = & (1+i)\{\frac{1}{4}(\zeta_{j_{\alpha}+1} + 2\zeta_{j} + \zeta_{j_{\alpha}-1})F^{h}(q_{l}, D^{1}q_{l})  -  \\
&  & \mbox{} 2(D_{+\alpha}\zeta_{j} + D_{-\alpha}\zeta_{j})(D_{+\alpha}q_{j} + D_{-\alpha}q_{j})-\frac{1}{4}(q_{j_{\alpha}+1} + 2q_{j} + q_{j_{\alpha}-1}^{h})\Delta^{h}\zeta_{j}\} ,
\end{eqnarray*}
Applying Duhamel. 
\begin{eqnarray*}
\lefteqn{\displaystyle \|q(t)\zeta)\|_{\Lh{4}(\R{2})}} \\
& \stackrel{\Delta}{=} &\sum_{k} \|\qplusk{}\zeta\|_{\Lh{4}(\R{2})}+ \|\qmink{}\zeta\|_{\Lh{4}(\R{2})} \\
&\letilde &  \|\mathcal{U}^h(t-\tau)(q(\tau)\zeta)\|_{\Lh{4}(\R{2})} + \int_{\tau}^{t} \frac{1}{(t-s)^{1/2}}\{\|(|q(s)|^{3}  +\\
& & \mbox{} +|D^{1}q(s)||q(s)|)\zeta\|_{\Lh{4/3}(\R{2})} + \||D^{1}q(s)||D^{1}\zeta|\|_{\Lh{4/3}(\R{2})}+  \\
& & \mbox{} + \||q||\Delta^{h}\zeta\||_{\Lh{4/3}(\R{2})} \} ds \\
&\letilde &  \|q(\tau)\zeta\|_{\Lh{4}(\R{2})} + \int_{\tau}^{t} \frac{1}{(t-s)^{1/2}}\{\|(|q(s)|^{3} +|D^{1}q(s)||q(s)|)\zeta\|_{\Lh{4/3}(\R{2})} \\
& & \mbox{}  + \||D^{1}q(s)||D^{1}\zeta|\|_{\Lh{4/3}(\R{2})}+  \||q||\Delta^{h}\zeta\||_{\Lh{4/3}(\R{2})} \} ds 
\end{eqnarray*}
By definition, $\displaystyle \qpluski{j}{a} = D_{+k}\uh{j}\cdot e_{j}^a$, so 
\[
D_{+i}\qpluski{j}{a}= O(|D^{2}\uh{j}|, |D^{1}u_{j}|^{2}).
\]
The $|D^{1}u_{j}|^{2}$ term can be seen as the error that appears from using $q$ to estimate
$D^{1}u^{h}$;  if it were not there, we would immediately have that $|D^{1}q| = O(|D^{2}u^{h}|) $ and then exploit the $L^{2}(\Lh{2})$ bounds on $u^{h}$ to attain bounds on $q\zeta$ from the above integral inequality.  Instead, we keep in mind that the use of $q$ was purely to make use of the linear estimates;  it is $D^{1}u^{h}\zeta$ that we need to bound. So long as $q$ satisfactorily controls $D^{1}u^{h}$,  we can achieve the desired bounds.  This control is a  consequence of Observation 4; namely,
\[
\begin{array}{ccc}
|D_{+k}\uh{j}|^{2} = O(|\qplusk{j}|^{2}), &  & |D_{-k}\uh{j}|^{2} = O(|\qmink{j}|^{2}).
\end{array}
\]
\begin{eqnarray*}
\lefteqn{\|D^{1}u^{h}\zeta\|_{\Lh{4}(\R{2})}}\\
& \letilde & \|q(t)\zeta)\|_{\Lh{4}(\R{2})}\\
&\letilde &  \|D^{1}u^{h}(\tau)\zeta\|_{\Lh{4}(\R{2})} + \int_{\tau}^{t} \frac{1}{(t-s)^{1/2}}\{ \|D^{1}u^{h}(s)\|_{\Lh{2}(\R{2})}\||D^{1}u^{h}(s)|^{2}\zeta\|_{\Lh{4}(\R{2})}  \\
& & \mbox{}  + \|(|D^{2}u^{h}(s)| + |D^{1}u^{h}(s)|^{2})q(s)\zeta\|_{\Lh{4/3}(\R{2})}  \\
& & \mbox{}  + \|(|D^{2}u^{h}(s)| + |D^{1}u^{h}(s)|^{2})D^{1}\zeta\|_{\Lh{4/3}(\R{2})}ds + C(R, E^h[f^h]) \\
&\letilde & \|D^{1}u^{h}(\tau)\zeta\|_{\Lh{4}(\R{2})}   + \int_{\tau}^{t} \frac{1}{(t-s)^{1/2}}\|D^{2}u^{h}(s)\|_{\Lh{2}(B_{2R}(x_{0}))}\|D^{1}u^{h}\zeta\|_{\Lh{4}(\R{2})}   \\
& & \mbox{}  + \|D^{2}u^{h}(s)\|_{\Lh{2}(B_{2R})(x_{0})}\|\|D^{1}u^{h}(s)\zeta\|_{\Lh{4}(\R{2})} \\
& & \mbox{}  + \frac{1}{R^{2}}\|D^{2}u^{h}(s)\|_{\Lh{2}(B_{2R})(x_{0})} + \||D^{1}u^{h}(s)|^{2}D^{1}\zeta|^{2}\|_{\Lh{2}(\R{2})} ds + C(R, E^h[f^h])\\
&\letilde & \|D^{1}u^{h}(\tau)\zeta\|_{\Lh{4}(\R{2})}   + \int_{\tau}^{t} \frac{1}{(t-s)^{1/2}}\{\|D^{2}u^{h}(s)\|_{\Lh{2}(B_{2R}(x_{0}))}\|D^{1}u^{h}(s)\zeta\|_{\Lh{4}(\R{2})} \\
& & \mbox{} +\frac{1}{R^{2}}\|D^{2}u^{h}(s)\|_{\Lh{2}(B_{2R})(x_{0})} \} ds + C(R, E^h[f^h])
\end{eqnarray*}  
Now taking $L^{r}$ in time of both sides, we know by  Young's Inequality that we can take $\|D^{2}u^{h}\|_{L_{h}^{2}(B_{2R}(x_{0}))}$ in $L^{2}$ (this term is bounded) and $\|D^{1}(u^{h}(s)\zeta)\|_{L_{h}^{4}}$ in $L^{r}$.  Integrating on a time interval of $\delta R^{2}$ and using the $L^{2}(H_{h}^{2})$ bounds, the second term is easily bounded.    Therefore, 
\[
\|D^{1}(u^{h}\zeta)\|_{L^{r}([\tau, t_{0}+\delta R^{2}], L_{h}^{4}(\R{2})} \leq  C(\|u^{h}(\tau)\|_{H_{h}^{2}}, E^{h}[f^{h}]).
\]
\end{subproof}
\item[(iv)] \fbox{$L^{\infty}(H_{h}^{2})$} 
\begin{subproof}
Since
\[
D^{1}\qpluski{j}{a} = O(|D^{2}\uh{j}|, |D^{1}u_{j}|^{2})
\]
and we have from (iii)  that $u_{h}$ has uniform $L^{r}(W_{h}^{1,4})$ bounds in a 
reduced parabolic cylinder, it suffices to attain uniform, local $L^{\infty}(W_{h}^{1,\infty})$ bounds on $q\zeta$.  
\begin{eqnarray*}
\lefteqn{\displaystyle \|D^{1}q(t)\zeta)\|_{\Lh{\infty}(\R{2})}} \\
&\letilde &  \|\mathcal{U}^h(t-\tau)(D^1q(\tau)\zeta)\|_{\Lh{\infty}(\R{2})} + \int_{\tau}^{t} \frac{1}{(t-s)^{3/4}}\{\|(|q(s)|^{3} + \\
& & \mbox{} +|D^{1}q(s)||q(s)|)\zeta\|_{\Lh{4/3}(\R{2})}  + \|D^{1}q(s)||D^{1}\zeta|\|_{\Lh{4/3}(\R{2})}\} ds + C(R, E^h[f^h])\\
&\letilde &  \|D^{1}q(\tau)\zeta\|_{\Lh{\infty}(\R{2})} + \int_{\tau}^{t} \frac{1}{(t-s)^{3/4}}\{\|q(s)\|_{\Lh{2}(\R{2})}\||q(s)|^{2}\zeta\|_{\Lh{4}(\R{2})} \\
& & \mbox{} +\|D^{1}(|q(s)|^{2}\zeta)\|_{\Lh{4/3}(\R{2})}  + \|(|D^{1}q(s)| + |q|)|D^{1}\zeta|\|_{\Lh{4/3}(\R{2})} \}+ \\
& & \mbox{} +C(R, E^h[f^h]) ds 
\end{eqnarray*}
\begin{eqnarray*}
&\letilde &  \|D^{1}q(\tau)\zeta\|_{\Lh{\infty}(\R{2})} + \int_{\tau}^{t} \frac{1}{(t-s)^{3/4}}\{\|q(s)\|_{\Lh{2}(\R{2})}\||q(s)|^{2}\zeta\|_{\Lh{4}(\R{2})}  + \\
& & \mbox{}  \||q(s)|^{2}\zeta\|_{\Lh{4}(\R{2})}  + \|(|D^{1}q(s)| + |q|)|D^{1}\zeta|\|_{\Lh{4/3}(\R{2})} \} ds \\ 
& & \mbox{} +C(R, E^h[f^h]) ds \\
&\letilde &  \|D^{1}q(\tau)\zeta\|_{\Lh{\infty}(\R{2})} + \int_{\tau}^{t} \frac{1}{(t-s)^{3/4}}
\|q(s)\|_{\Lh{4}(B_{2R}(x_{0})}\|D^{1}q(s)\zeta\|_{\Lh{2}(\R{2})}ds + \\
& & \mbox{  } +  C(R, E^h[f^h]) 
\end{eqnarray*}
Taking $L^{\infty}$ in time of both sides, we have by Young's Inequality that we can take
$ \|q(s)\|_{\Lh{4}(B_{2R}(x_{0})}$ in $L^{r}$ for any $r>4$ to achieve the bound
\[
\|D^{1}(u^{h}\zeta)\|_{L^{\infty}([\tau, t_{0}+\delta R^{2}], L_{h}^{\infty}(\R{2})} \leq  C(\|u^{h}(\tau)\|_{H_{h}^{2}}, E^{h}[f^{h}]).
\]
\end{subproof}
\end{itemize}
\end{proof}

\begin{proof} (\emph{Theorem \ref{fellgthm1}})
\subsubsection*{Step 1: Attaining Bounds}
This is the content of Lemma $\ref{fellgthm2}$.
\subsubsection*{Step 2:  Estimating the size of the concentration set}
There is no change in this argument from the case of the harmonic map heat flow.
\subsubsection*{Step 3:  Convergence of the Interpolants}
The argument has been illustrated  in the case of the harmonic map heat flow.  The only work that
has to be done  to accommodate for the presence of the Schr\"odinger term is to show that on a compact set $Q$ for which we have uniform $L^{\infty}(H_{h}^2)$ bounds, that the interpolant $p_h$ converges.  For $\phi \in C_{0}^{\infty}(Q)$, we need to show that
\[
\int_{Q} \partial_{t} p^h \cdot \phi + \nabla p^h \cdot \nabla \phi - \nabla p^h \cdot \nabla \nu (p^h) \phi+ \sum_{i} \{(\nu(p^h) \wedge \partial_{i} p^h) \cdot \partial_{i} \phi - (\partial_{i} \nu (p^h) \wedge \partial_{i} p^h)\cdot \phi\} = O(h)
\]
Expanding the $\partial_{t}p^h \cdot \phi$ term,  collecting terms in $\partial_{t} u^h$ and using
the $L^2(L_{h}^2)$ bound on $\partial_t u^h$, 
the only term that is not already $O(h)$ is $\partial_t \uh{j} \cdot \phi$ which contains the following terms:
\[
\partial_{t}\uh{j}(t)\cdot \phi =  (\uh{j}\wedge \Delta \uh{j}) \cdot \phi +  (\Delta^{h}\uh{j}(t) + \lambdah{j} \nu_{j})\cdot \phi.
\]

It was shown in the harmonic map heat flow case that the terms involving the damping term have the correct cancellation. 
It remains to show that
\begin{eqnarray*}
\lefteqn{\int_{Q}  (\uh{j} \wedge \Delta^h \uh{j}) \cdot \phi+ \sum_{i} \{(\nu(p^h) \wedge \partial_{i} p^h) \cdot \partial_{i} \phi - (\partial_{i} \nu (p^h) \wedge \partial_{i} p^h)\cdot \phi \}} \\
& = & \int_{Q} \sum_{i} \{(D_{+i}(\nu_{j}  \wedge D_{-i} \uh{j}) - D_{+i} \nu_{j} \wedge D_{+i}\uh{j})\cdot \phi + \\
&  & \mbox{ }(\nu(p^h) \wedge \partial_{i} p^h) \cdot \partial_{i} \phi - (\partial_{i} \nu (p^h) \wedge \partial_{i} p^h)\cdot \phi \}\\
& = &O(h)
\end{eqnarray*}
The uniform $L^{\infty}(H_{h}^{2})$ bounds resolves most of the terms.   In addition, the following manipulation
\begin{eqnarray*}
\lefteqn{\int _Q  \sum_i D_{+i}(\nu_{j}  \wedge D_{-i} \uh{j}) \cdot \phi + (\nu(p^h) \wedge \partial_{i} p^h) \cdot \partial_{i} \phi} \\
& = & \int _Q  \sum_i((\nu_j - \nu(p^h))\wedge D_{-i}u_j)\cdot \phi  - \nu(p^h)\wedge (h \Delta^h \uh{j}) \cdot \phi  + O(h) 
 \end{eqnarray*}
combined with these bounds permits us to conclude that $p_h$ converges in $Q$.   The argument to show that $p$ solves the equation on $\Rplus \times \R{2}$ proceeds with no change from the harmonic map heat flow case. 
\end{proof}

A consequence of the proof for Theorem \ref{fellgthm1} is the special case for small energy initial data.
\begin{corollary}
\label{smallllg}
There is a constant $\epsilon_{0}$ such that for $f \in H^{1}(\R{2}, \mathcal{N})$, $E[f] < \epsilon_0$, there exists a smooth global solution to  $(\ref{llgequation})$.  
\end{corollary}

\vspace{0.3in}
\noindent \textbf{Acknowledgements}  This work was part of a dissertation completed at the Courant Institute.  I wish to thank  Fanghua Lin for suggesting a problem that led to this work and Jalal Shatah for his interest and guidance throughout the process.

\section*{\center{\Large{Appendix A}}\markboth{Appendix}{Appendix}}
\setcounter{equation}{0}
\renewcommand{\theequation}{A.\arabic{equation}}


Let $\Omega \subset \R{2}$ and  $\zeta \in C_{0}^{\infty}(\Omega)$. Then
\begin{equation}
\label{clocsob-int}
\||f|^{2}\zeta\|_{L^{2s}(\R{2})} \leq C(s)\|f\|_{L^{p}(\Omega)} \|\nabla(f\zeta)\|_{L^{q}(\R{2})}; \  \
   \frac{s+1}{2s} = \frac{1}{p} + \frac{1}{q}.
\end{equation}
\begin{proof} We first take the easy case $q < 2$.  Then by H\"older and Sobolev embedding, 
\[
\||f| |f \zeta|\|_{L^{2s}(\Omega)} \leq C\|f\|_{L^p(\Omega)}\|f \zeta\|_{L^r(\R{2})} \leq C\|f\|_{L^p(\Omega)}\|\nabla (f\zeta)\|_{L^q(\R{2})}
\]  
the last inequality holding iff $\frac{1}{q} - \frac{1}{2} > 0$ iff $q <2$.  $p, q$ and $s$ are related as follows:
\[
\frac{1}{p} = \frac{1}{2s} - \frac{1}{r} = \frac{1}{2s} - \frac{1}{q} + \frac{1}{2}
\]
This is the inequality that we want.   

Now consider the case $q > 2$. The proof follows in the same manner as that
for the normal Gagliardo-Sobolev-Nirenberg inequalities, with the help of the following identity:
\[
\partial_{k}(f^{2s}\zeta^{s}) = 2s\partial_{k}(f\zeta)f^{2s-1}\zeta^{s-1}-f^{2s}\partial_{k}(\zeta^{s}).
\]
We can write
\begin{eqnarray*}
\lefteqn{|\partial_{1}(f^{2s}(x_1, x_2)\zeta^{s}(x_1,x_2))|}\\
& \leq & 2s 
\int_{-\infty}^{x_{1}}|\partial_{1}(f(y_1,x_2)\zeta(y_1,x_2))f^{2s-1}(y_1,x_2)\zeta^{s-1}(y_1,x_2)| dy_1 +\\
&  &\mbox{ } 
\int_{-\infty}^{x_{1}}|f^{2s}(y_1,x_2)\partial_{1}(\zeta^{s}(y_1,x_2))|dy_1
\end{eqnarray*}
and a similar expression holding $x_{1}$ fixed and integrating on the interval $(-\infty, x_{2})$.  
Multiplying these two expressions together,  and integrating, we have
\begin{eqnarray*}
\lefteqn{\||f|^{2}\zeta\|_{L^{2s}}^{2s} =  \int_{-\infty}^{\infty} \int_{-\infty}^{\infty} |f^{2s}(x_1,x_2)\zeta^s(x_1,x_2)|^{2}dx_1 dx_2}\\
& \leq & C(s)\{(\int_{-\infty}^{\infty} \int_{-\infty}^{\infty} |\nabla(f\zeta)||f^{2s-1}||\zeta^{s-1}|dx_1 dx_2)^{2} + \\ & & \mbox { } +(\int_{-\infty}^{\infty} \int_{-\infty}^{\infty} |\nabla(f\zeta)||f^{2s-1}||\zeta^{s-1}|dx_1 dx_2)(\int_{-\infty}^{\infty} \int_{-\infty}^{\infty} \nabla|\zeta^{s}||f^{2s}|dx_1 dx_2) + \\
& & \mbox{ } +(\int_{-\infty}^{\infty} \int_{-\infty}^{\infty} |\nabla(\zeta^s)||f^{2s}| dx_1 dx_2)^{2}\} \\
& \leq & C(s) \|\nabla(f\zeta)\|_{L^{q}}^{2}\|f\|_{L^{p}}^{2}\||f|^{2}\zeta\|_{L^{2s}}^{s-1}
\end{eqnarray*}
\end{proof}

\section*{\center{\Large{Appendix B}}\markboth{Appendix}{Appendix}}
 Let $\mathcal{N}$ be a compact smooth manifold.  Given  $p \in \mathcal{N}$, there exists a neighborhood $U$ of $p$ and a local coordinate system $\{y^i\}$ on $U$.  An immediate basis for the tangent space at any point in $U$ is given by the coordinate vector fields $X_i = \fracpartial{}{y^i}$.   However, it is sometimes not convenient to use vectors that come 
from a coordinate system.   For any $n$ linearly independent vectors in $\TyN$,  we denote the value of the field in $p$ by
\[
X(p) = \{X_1(p), \ldots, X_n (p)\}.
\]
$X(p)$ is a basis for $\TpN$.  
A  \emph{moving frame} on $\mathcal{N}$ is a  function $F:  p \in \mathcal{N} \rightarrow X(p)$, where  $F(p)$ defines a basis for $\TpN$.   If $F$ is defined for all $p \in \mathcal{N}$ then $F$ is a \emph{global frame} for $\mathcal{N}$.

Throughout the process of linearizing the discrete LLG in Section $\ref{linearizedsystem}$ , we assumed the existence of a smooth global frame on $\mathcal{N}$.  The problem is that this assumption can
be thwarted by simple examples,  the sphere $S^2$ being one such case.  Fortunately, this potentially unfortunate hindrance can be successfully bypassed by at least  two methods.  These methods are standard  so we provide here a highly customized sketch of these two remedies:  (1) embedding $\mathcal{N}$ isometrically as a totally geodesic submanifold into a larger space which does admit a global frame and (2) using a pullback frame on $\R{2} \times \R{}$ by exploiting the smoothness of the map $u$.  A good reference for this is \cite{helein}. 
\begin{itemize}
\item[(1)] \noindent \emph{Embedding into something larger}

\noindent One way of introducing a global frame is to consider the equations on a tubular neighborhood $U$ of $\mathcal{N}$ in $\R{3}$ which does admit such a frame.  In this case $y \in U$  can be represented by coordinates $(u(y), \rho(y))$.  For  $u \in \mathcal{N}$,  $\nu(u)$ can be extended to $U$ by $\nu(y) = \nu(u(y))$.  Then for small $\rho$, 
\[
y = u + \rho \nu(y). 
\] 
The equation for $y$ is given by 
\[
\partial_t y = \nu(y) \wedge \Delta y +\alpha \nu(y) \wedge (\nu(y) \wedge \Delta y).
\]
All derivatives of $u$ can be computed using a global frame on $\R{3}$.  
However, this procedure has the drawback of losing all the geometric structure of the problem and
consequently all the estimates and cancellations these structures imply.  

We can view the LLG as a specific example of a \emph{Schr\"odinger map}  which is a map $u: \R{d} \times \R{} \rightarrow (\mathcal{N}, g, J)$, where $\mathcal{N}$ is a K\"ahler manifold with a complex
structure $J = \nu \times$ and a metric $g$ giving rise to covariant differentiation $D$, satisfying
\[
\partial_t u = J(u) D\partial u + \alpha D\partial u.
\]
What we really need  is to embed $\mathcal{N}$ isometrically into a larger manifold $\tilde{\mathcal{N}}$  which itself has a complex structure $\tilde{J}$ where $\tilde{J}|_{\mathcal{N}} = J$ and in such a way that the equation in $\tilde{\mathcal{N}}$ 
\[
\partial_t \tilde u = \tilde{J}(\tilde{u}) D\partial \tilde{u} + \alpha D\partial \tilde{u}
\]
has $\mathcal{N}$ as an invariant subspace.
For our problem this can be achieved in the following manner.
Figure \ref{embedding} attempts to illustrate this embedding.  First we embed $\mathcal{N} \subset \R{3} \subset \R{4}$ (a complex manifold must be of even dimension.)  
\begin{figure}[htbp]
\begin{center}
\includegraphics[height =.7\textwidth, angle = 270]{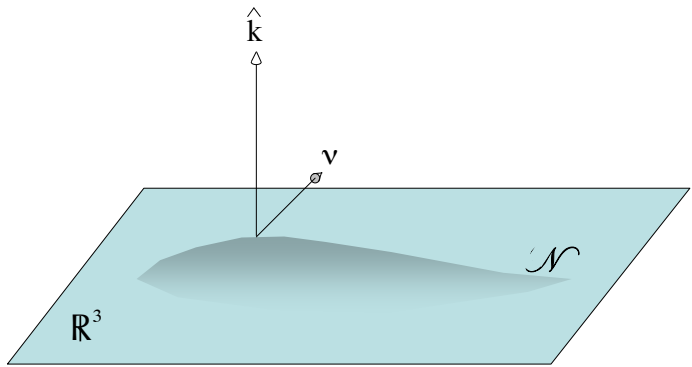} 
\caption{Embedding $\mathcal{N} \subset \R{3} \subset \R{4}$ } 
\label{embedding}
\end{center}
\end{figure}
On a tubular neighborhood $U$ of $\mathcal{N}$  in $\R{4}$ we introduce the metric $\tilde{g}$ given
in coordinates $\tilde{u} = u + \rho \nu(u) + \xi \hat{k}$, $u\in \mathcal{N}$:
\[
ds^2 = ds_{\mathcal{N}}^2 + d\rho^2 + d\xi^2
\] 
and the almost complex structure $\tilde{J}(\tilde{u})$ 
\[
\tilde{J}\tilde{y} = \tilde{J}(y(u) + a \nu(u) + b\hat{k}) = Jy(u) - b\nu(u) + a\hat{k}.
\]
The four dimensional manifold $U$  has a global frame  by construction and is given by a metric  giving rising to a covariant derivative $\tilde{D}$ such that
\[
\tilde{D} \tilde{J} = 0. 
\]
The equation
\[
\partial_t y = \tilde{J}\tilde{D}\partial y+ \alpha \tilde{D} \partial y
\]
has $\mathcal{N}$ as an invariant subspace.  This last statement is true since for  $u_0 \in \mathcal{N}$, 
we have $\rho_0 = \xi_0$,  and the equation  reduces to 
\begin{eqnarray*}
\left\{\begin{array}{lcl}
\displaystyle
\partial_t u = J(u)D\partial u+ \alpha D \partial u \\
u( 0) =  u_0
\end{array}\right.
\end{eqnarray*}

\item[(2)] \noindent \emph{Pullback Frame} 

\noindent Consider a smooth embedded surface $\mN \subset \R{3}$.  A pullback frame for a map $u: \Omega \subset \R{2} \rightarrow \mN$ is a choice, for each point $x\in \Omega$,  of an orthonormal basis $\{e_1(x), e_2(x)\}$ for $\mathcal{T}_{u(x)}\mN$. Since $\mN$ is embedded, there is a well-defined tangent plane $\mathcal{T}_{u(x)}\mN$ for every point $x \in \Omega$.  An orthonormal basis can always be specified at a point so this type of frame  always exists.   If $\mathcal{N}$ admits a global frame, then a pullback frame for $u$ automatically exists, given by $\{e_1(u(x)), e_2 (u(x))\}$.   However, the converse is not true.  The pullback frame is a frame on the domain, and might be different at $u(x_1)$ and $u(x_2)$ even when the two values coincide on $\mathcal{N}$;   a pullback frame, therefore, might exist in the absence of a global frame.  

The construction of such a frame is most apparent in the case $u: I \subset \R{} \rightarrow \mathcal{N}$.  
Let $e_0$ be
a frame at $T_{u(0)} \mN$.  Then we can solve  for $e$ by parallel transport (i.e., $e$ satisfying $D_x e = 0$) which gives the following o.d.e.: 
\begin{equation}
\label{parallel transport}
\left\{\begin{array}{lcl}
\displaystyle
\partial_x e + (e \cdot (\nu(u))_x) \nu(x) = 0 \\
e( 0) =  e_0
\end{array}\right.
\end{equation}
Figure \ref{pt} illustrates this selection of a frame.
\begin{figure}[htbp]
\begin{center}
\includegraphics[width =.75\textwidth]{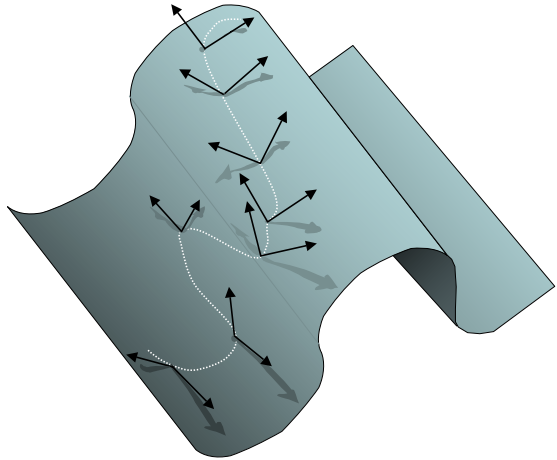} 
\caption{Constructing $\{e_1, e_2\}$ by parallel transport} 
\label{pt}
\end{center}
\end{figure}

From this construction it is apparent that
the frame reflects the geometry of $u$, since $e$ and $u$ have the same regularity. 
Higher dimensional analogs of this construction are standard.
\end{itemize}

\bibliography{bibllg}
\bibliographystyle{plain}

\end{document}